\documentclass[10pt]{amsart}
\usepackage{amssymb}
\usepackage{MnSymbol}

\usepackage{amscd}
\usepackage{amsthm}
\usepackage[dvips]{graphicx}
\usepackage{color}
\usepackage[all]{xy}
\usepackage{verbatim}

\newtheorem{Lemma1}{{Lemma}}[section]
\newtheorem{Theo1}[Lemma1]{{Theorem}}
\newtheorem*{Theo2}{{Theorem}}
\newtheorem{Def1}[Lemma1]{{Definition}}
\newtheorem{Prop1}[Lemma1]{{Proposition}}
\newtheorem{Claim1}[Lemma1]{{Claim}}
\newtheorem{Rem1}[Lemma1]{{Remark}}
\newtheorem{Cor1}[Lemma1]{{Corollary}}
\newtheorem{Ex1}[Lemma1]{{Example}}
\newtheorem{Not1}[Lemma1]{{Notation}}

\newenvironment{Lemma}{\begin{Lemma1}}{\end{Lemma1}}
\newenvironment{Def}{\begin{Def1}\rm}{\end{Def1}}
\newenvironment{Prop}{\begin{Prop1}}{\end{Prop1}}
\newenvironment{Rem}{\begin{Rem1}\rm}{\end{Rem1}}
\newenvironment{Theorem}{\begin{Theo1}}{\end{Theo1}}

\newenvironment{Cor}{\begin{Cor1}}{\end{Cor1}}

\newenvironment{Example}{\begin{Ex1}\rm}{\end{Ex1}}

\newcommand{\red}{ }

\title[Ore localisation and Goldie's theorem of dg-algebras]{Ore Localisation for differential graded rings; Towards
Goldie's theorem for differential graded algebras}

\author{Alexander Zimmermann}
\address{\newline
Universit\'e de Picardie,
\newline D\'epartement de Math\'ematiques et LAMFA (UMR 7352 du CNRS),
\newline 33 rue St Leu,
\newline F-80039 Amiens Cedex 1,
\newline France}
\email{alexander.zimmermann@u-picardie.fr}

\date{November 28, 2023}

\newcommand{\dar}{\downarrow}
\newcommand{\lra}{\longrightarrow}

\newcommand{\ra}{\rightarrow}
\newcommand{\sdp}{\times\kern-.2em\vrule height1.1ex depth-.05ex}
\newcommand{\epi}{\lra \kern-.8em\ra}

\newcommand{\N}{{\mathbb N}}
\newcommand{\R}{{\mathbb R}}

\newcommand{\ol}{\overline}

\newcommand{\Z}{{\mathbb Z}}

\newcommand{\soc}{\textup{soc}}

\newcommand{\im}{\textup{im}}
\newcommand{\dickebox}{{\vrule height5pt width5pt depth0pt}}

\newcommand{\rudim}{\textup{rudim}}
\newcommand{\ludim}{\textup{ludim}}

\newcommand{\dgnil}{\textup{dgnil}}

\newcommand{\Prad}{\textup{Prad}}
\newcommand{\Nil}{\textup{Nil}}
\newcommand{\dgudim}{\textup{dg-udim}}
\newcommand{\lann}{\textup{lann}}
\newcommand{\rann}{\textup{rann}}
\newcommand{\ann}{\textup{ann}}
\newcommand{\End}{\textup{End}}

 \normalbaselines

\begin{document}

\begin{abstract}
We study Ore localisation of differential graded algebras. Further we 
define dg-prime rings, dg-semiprime rings, and study the dg-nil radical of dg-rings. Then, we 
define dg-essential submodules, dg-uniform dimension, and apply all this 
to a dg-version of Goldie's theorem on prime dg-rings. 
\end{abstract}

\maketitle

\section*{Introduction}

Differential graded algebras first appeared  in a paper by Cartan~\cite{Cartandg} 
and were then developed mainly in the context of algebraic topology, algebraic geometry and 
differential geometry. For a modern treatment we refer to Yekutieli~\cite{Yekutieli}. 
For a commutative base ring $K$ a 
differential graded $K$-algebra is an associative unital $\Z$-graded $K$-algebra
together with a degree $1$ endomorphism $d$ of square $0$ satisfying 
Leibniz formula $$d(a\cdot b)=d(a)\cdot b+(-1)^{|a|}a\cdot d(b)$$ 
for all $a,b\in A$ such that $a$ is homogeneous of degree $|a|$.  
So a differential graded algebra is an algebra at first, and this 
fact remained largely unexplored until very recently. 
In a recent sequel of papers by the author \cite{dgorders, dgBrauer}, by Orlov 
\cite{Orlov1,Orlov2}, by Aldrich and 
Garcia-Rozas~\cite{Tempest-Garcia-Rochas}, and by Goodbody~\cite{Goodbody}
differential graded rings are considered from a ring theoretic point of view.
Still Orlov's papers have an algebraic geometric perspective in mind. 

Orlov  \cite{Orlov1,Orlov2}, and later but independently  \cite{dgorders}, gave 
a definition of semisimple dg-algebras 
from an algebra point of view, namely that in a dg-simple dg-algebra there is no non-trivial 
two-sided differential graded ideal, and a dg-semisimple dg-algebra is a direct product of 
dg-simple dg-algebras. In contrast to the classical case, this 
definition leads to a 
concept different from when we consider dg-algebras whose dg-module category 
is semisimple. This latter point of view has been completely settled by Aldrich and 
Garcia-Rozas~\cite{Tempest-Garcia-Rochas}. Following Orlov's approach, a definition of a dg-radical is 
given by Goodbody~\cite{Goodbody}, including a dg-version of Nakayama's lemma.
In \cite{dgorders} there is a different, within many respects
more powerful study of a dg-radical and Nakayama's lemma, and closes with a definition 
of dg-orders and dg-class groups including many properties. In \cite{dgBrauer} 
a dg-version of the Brauer group is given as well as the proof that the dg-Brauer group
coincides with the classical one. 

The present paper continues these investigations. 

In classical non-commutative ring theory Ore localisation is a very important tool.
Using it, an important result is Goldie's theorem for prime rings satisfying the ascending
chain condition on right annihilator ideals and satisfying that there is only a finite
number of uniform ideals whose sum is direct and forms an essential ideal of
the ring. This number is well-defined
and is called the uniform dimension of the ring. Rings that satisfy these conditions
are called  Goldie rings. Goldie's theorem \cite{Goldie, Goldie2}
states that the Ore localisation at the regular elements
of a prime Goldie rings is a simple Artinian ring.

A group-graded version for an abelian group
of this theorem was given in its most accessible form by Goodearl
and Stafford~\cite{GoodearlStafford}.

In the present note we first consider Ore localisation and show that the Ore
localisation at homogeneous regular elements of a differential graded algebra still
is a differential graded algebra, extending the differential graded structure of the initial algebra. 
A different approach of Ore localisation
focusing on the homology of the dg-algebra was given by Braun, Chuang and Lazarev~\cite{BraunChuangLazarev}.
The authors of \cite{BraunChuangLazarev} mention that their construction is difficult to realise explicitly. 
Our construction provides such an explicit construction (cf Proposition~\ref{BCLisquasiiso})
but also generalises the construction considerably. 

Along the way we consider the dg-nil radical of an algebra and see that the dg-prime radical, which is the
intersection of dg-prime ideals, in general strictly includes the dg-nil radical. {\red Just as in the classical situation we show that if $(A,d)$ satisfies the artinian condition on twosided differential graded ideals, then the dg-nil radical coincides with the dg-Jacobson radical in its twosided version and with the dg-prime radical.}

Furthermore, we define naturally and study elementary properties of dg-essential ideals, and the
naturally occurring dg-uniform dimension. In the classical theory 
the dg-right singular ideal is the set of left annihilators of
essential right ideals. In contrast to the classical or the graded case the differential graded version 
is not a two-sided ideal in
general, nor is it nilpotent. We study the relation with the classical and with the graded
case, in particular related to these concepts for the homology algebra and the cycles.

Finally, we prove a dg-version of Goldie's theorem for a dg-algebra $(A,d)$,
basically under the hypotheses that the cycles $\ker(d)$ satisfy the hypotheses of
the graded version \cite{GoodearlStafford} of Goldie's theorem, where simplicity is defined 
as the absence of non-trivial two-sided dg-ideals. We show that this 
hypothesis is stronger than asking the dg-version of the corresponding
Goldie's hypotheses. However, since the dg-singular ideal is badly behaved, it is unlikely
that a direct generalisation of Goldie's theorem to the dg-world is possible.

The paper is structured as follows. In Section~\ref{recall} we recall the necessary 
definitions and concepts of the classical and of the graded theory around Goldie's theorem.
In Section~\ref{orerevisited} we revise existing results around Ore localisation. Apart from the 
classical definition and existence result, we review Braun, Chuang, Lazarev's result 
\cite{BraunChuangLazarev} on  Ore localisation of the homology of a dg-algebra.
We also  prove our first main result Theorem~\ref{orelocalisationofdg}, which states
that the differential of a dg-algebra
extends to the structure of a dg-algebra on the Ore localisation of a dg-algebra 
at homogeneous regular elements.  We then draw immediate consequences on the classical 
Goldie theorem and study consequences for the two available versions of localisation, namely
the version from \cite{BraunChuangLazarev} and the version from Theorem~\ref{orelocalisationofdg}.
Section~\ref{dgprimesect} then studies questions on a version of prime and semiprime dg-algebras
in the differential graded sense. In connection with this we study properties of the dg-nil radical and discover
that in contrast to the classical case, the prime radical, i.e. the intersection of dg-prime ideals 
strictly contains the maximal nilpotent dg-ideal. In Section~\ref{dgessential} we define and study
elementary properties of a dg-version of a module that is essential as dg-module.  Using the results of
Section~\ref{dgessential} we develop in Section~\ref{dguniformsect} the theory of dg-uniform dimension, 
show that it is well-defined and define and study the dg-singular ideal. We show by an example
that this dg-singular ideal is not two-sided, unlike the graded case, or the classical case, and 
give links to the singular ideal of the homology, respectively its relation to cycles and the 
singular ideal of the cycles. Section~\ref{dgleftrightannihilsect} we consider differential graded 
left or right annihilator ideals and study the connections with the dg-singular ideal and dg-essential 
ideals. Finally, Section~\ref{dgGoldietheoremsect} we prove our second main result 
Theorem~\ref{essentialcontainsregular}, namely a 
dg-version of Goldie's theorem.

\section{The classical situation: Goldie's theorem and its graded version}

\label{recall}

\subsection{The classical ungraded case}

We refer to McConnell and Robson~\cite[Chapter 2]{McconnellRobson} for the treatment of the
present section.

We start to recall a number of standard definitions in ring theory.

\begin{itemize}
\item
A ring is called {\em semiprime} if it does not have any nilpotent right ideal, which is actually equivalent with
$R$ not having any two-sided nilpotent ideal.
\item
A ring $R$ is {\em prime} if for any two non zero two-sided ideals
$I$ and $J$ of $R$ one has $IJ\neq 0$.
\item
A $R$ right ideal $I$ is called {\em essential} if for any non zero $R$ right ideal $X$ one has
$I\cap X\neq 0$.
\item
The {\em right singular ideal $\zeta(R)$} is defined the as the set of $a\in R$ such that there is an
essential right ideal $E$ of $R$ with $aE=0$. This is an ideal indeed.
\item
Similarly, one defines an {\em essential submodule} of a module.
\item
An $R$-module $U$ is called {\em uniform} if it is non zero and any non zero submodule of $U$ is essential.
\item
A right module $M$ has {\em finite uniform dimension} if it contains no infinite direct
sum of non zero submodules.
\item
A right module $M$ with finite uniform dimension has {\em right uniform dimension $n$}
if any direct sum of submodules of $M$ has at most $n$ summands, and if there
is a direct sum of $n$ uniform submodules. This direct sum is essential in $M$.
Moreover, any direct sum of uniform submodules is essential if and only if it has $n$ summands.
\item
The {\em right uniform dimension $\rudim(R)$} of $R$ is defined to be the uniform dimension
of the right regular module $R$.
Similarly one defines the left uniform dimension  $\ludim(R)$.
\item A ring $R$ is a {\em right Goldie ring} if it has finite right uniform dimension
on right ideals and it satisfies the ascending chain condition of right annihilators.
\end{itemize}

With these preparations Goldie's celebrated theorem  states as follows
(cf e.g. \cite[2.3.6]{McconnellRobson}).

\begin{Theorem} (Goldie) \label{goldiestheoremclassical}
The following are equivalent.
\begin{itemize}
\item $R$ is semiprime right Goldie.
\item $R$ is semiprime, $\zeta(R)=0$ and the right uniform dimension of $R$ is finite.
\item The Ore localisation $Q$ of $R$ at the regular elements ${\mathcal C}(0)=:S$
of $R$ is semisimple artinian.
\end{itemize}
Further, in the last item $R$ is prime if and only if $Q$ is simple.
\end{Theorem}

%
%

\subsection{The graded situation}

Due to its crucial importance, Goldie's theorem was generalised in several directions.
We concentrate on graded rings for the moment. As a reference we use \cite{NastasescuVanOystaen}.

\begin{itemize}
\item
For a group $G$ we call a $G$-graded ring $R$ is a ring $R$ such that $R=\bigoplus_{g\in G}R_g$
such that $R_g\cdot R_h\subseteq R_{gh}$ for all $g,h\in G$. The ring is strongly graded if we
have $R_g\cdot R_h= R_{gh}$ for all $g,h\in G$.
\item
Let $R$ be a $G$-graded ring $R$. Then a $G$-graded module $M$ is an $R$-module
such that $M=\bigoplus_{g\in G}M_g$ and such that $R_g\cdot M_h\subseteq M_{gh}$ for all $g,h\in G$.
A $G$-graded module is gr-simple if it does not have any $G$-graded submodule.
It is gr-Artinian if any descending chain of $G$-graded submodules is finite.
\item
A $G$-graded ring is called gr-prime if for any two non zero $G$-graded two-sided ideals
$I$ and $J$ we have $IJ\neq 0$.
\item
For a $G$-graded module $M$ a $G$-graded submodule $N$ is gr-essential
if for any non zero $G$-graded submodule $X$ of $M$ we have $N\cap X\neq 0$.
\item
Likewise a $G$-graded module  is called gr-uniform if it is non zero
and any $G$-graded submodule is gr-essential.
\item
A $G$-graded module $M$ has finite gr-uniform dimension if it does not contain an
infinite direct sum
of $G$-graded non zero submodules.
\item
A $G$-graded ring is a right gr-Goldie ring if it has finite
right gr-uniform dimension on $G$-graded right ideals and it satisfies the ascending
chain condition on right annihilators of homogeneous elements.
\item The gr-singular ideal of $R$ is the set of homogeneous elements of $R$ such that there is
gr-essential right ideal $I$ with $aI=0$.
\end{itemize}

We should note that a prime ring is gr-prime, an essential submodule is
gr-essential submodule, etc.

\medskip

In this setting we recall from \cite[8.4.5 Theorem]{NastasescuVanOystaen}
the following result of
Goodearl and Stafford.

\begin{Theorem} \cite[Lemma 2, Theorem 1]{GoodearlStafford} \label{GoodearStafford}
Let $R$ be a $G$-graded ring, where $G$ is an abelian group. Suppose that
$R$ is a gr-prime right gr-Goldie ring.
\begin{enumerate}
\item
Then any non zero graded two-sided ideal of $R$ contains a non-nilpotent homogeneous element.
\item Further, the right gr-singular ideal of $R$ is nilpotent, whence $0$.
\item Finally the localisation of $R$
at homogeneous regular elements is a gr-simple, gr-Artinian ring.
\end{enumerate}
\end{Theorem}

Note that this is an almost faithful transposal of Goldie's theorem to the
graded situation, grading by an abelian group. According to the remarks
in the introduction of \cite{GoodearlStafford}, prior to this result
it was considered that Goldie's theorem does not transpose to the graded
situation, but this point of view is due to an misconception of the
result to be formulated.

\section{Ore localisation and differential graded rings}

\label{orerevisited}

\subsection{Classical Ore localisation revisited}

We recall the theory of Ore localisation from classical ring theory.

Let $S$ be a non empty multiplicatively closed subset of $R$ and define
$$\textup{ass}(S):=\{r\in R\;|\;\exists s\in S\::\;rs=0\}. $$
A {\em right quotient ring of $R$ with respect to $S$} is a ring $Q$
together with a ring homomorphism $\theta:R\lra Q$ such that
\begin{enumerate}
\item $\theta(S)\subseteq Q^\times$, the group of invertible elements in $Q$.
\item $\forall q\in Q\exists s\in S\exists r\in R\;:\;q\cdot\theta(s)=\theta(r)$
\item $\ker(\theta)=\textup{ass}(S)$.
\end{enumerate}
Similarly one defines the left quotient ring by modifying the second condition accordingly.

A multiplicatively closed subset $S$ satisfies the {\em right Ore condition} if
$$\forall r\in R\forall s\in S\exists r'\in R\exists s'\in S\;:\;rs'=sr'$$
Dually one defines the left Ore condition.
It is easy to see that if a right quotient ring exists, then $S$ satisfies the right Ore
condition. Further, by e.g. \cite[Chapter 2.1.12]{McconnellRobson}, if
the mutiplicatively closed set $S$ satisfies the right Ore condition, then $\textup{ass}(S)$ is a
two-sided ideal of $R$ and the right quotient ring $R_S$ with respect to $S$ exists if and only if
the image of $S$ in $R/\textup{ass}(S)$ consists of regular elements.

\subsection{Localisation of dg-algebras making homology classes invertible}

We recall results from \cite{BraunChuangLazarev} on localisation of
dg $k$-algebras and their dg-modules. The authors of \cite{BraunChuangLazarev} show that
for any set $T$ of regular elements of the homology of $R$ there is a dg-algebra
inverting the classes $T$.

For a subset $S$ of homogeneous elements of $H(A)$ the authors  \cite[Definition 3.1]{BraunChuangLazarev}
say that
a dg $A$-algebra $f:A\lra Y$ is $S$-inverting if for all $s\in S$
the homology class $f_*(s)\in H(Y)$ is invertible in the algebra $H(Y)$.
The derived localisation $L_S^{dgAlg}(A)$ is (cf \cite[Definition 3.3]{BraunChuangLazarev})
the initial object in the
full subcategory  on $S$-inverting dg algebras of $A\dar^{\mathbb L} dgAlg$.
Here the category  $A\dar^{\mathbb L} dgAlg$ is formed by dg algebra morphisms
$A\lra Y$ and morphisms being given by dg algebra homomorphism making commutative
the obvious commutative triangles.
Finally, the authors show \cite[Theorem 3.10]{BraunChuangLazarev} that
such a localisation exists and
$L_S^{dgAlg}(A)\simeq A\ast_{k(S)}^{\mathbb L}k(S,S^{-1}). $

The authors further show (cf remarks after \cite[Corollary 3.24]{BraunChuangLazarev})
that if $S$ is an Ore set in $H(A)$, then the homology of the localisation is the
localisation of the homology.

\subsection{Extending the differential to the Ore localisation}

In case $S$ is a multiplicative system with only regular elements in a differential graded ring $(R,d)$,
then the right quotient ring $R_S$ contains $R$, in the sense that $\theta:R\lra R_S$ is injective.
Hence any $q\in R_S$ can be written as $q=\theta(r)\theta(s)^{-1}$, or equivalently
$q\theta(s)=\theta(r)$. Hence, if $\widehat d$ is an extension of $d$ to $R_S$, and
if $q$ is a homogeneous element in $R_S$,
then
$$\theta(d(r))=\widehat d(\theta(r))=\widehat d(q\theta(s))=
\widehat d(q)\theta(s)+(-1)^{|q|}q\widehat d(\theta(s))=
\widehat d(q)\theta(s)+(-1)^{|q|}q\cdot \theta(d(s))$$
If $q\theta(s)=\theta(r)$ and $r$ is of degree $|r|$, and $s$
is of degree $|s|$, then $|q|=|r|-|s|$. Therefore
\begin{eqnarray*}\widehat d(q)&=&\theta(d(r))\cdot \theta(s)^{-1}-(-1)^{|q|}q\theta(d(s))\\
&=&\theta(d(r))\cdot \theta(s)^{-1}-(-1)^{|r|-|s|}\theta(r)\cdot \theta(s)^{-1}\cdot \theta(d(s))
\end{eqnarray*}
This formula holds also in case $\theta$ is not injective.
Hence, $\widehat d$, if it exists, is uniquely defined by the above formula.

We still need to show that this formula is well-defined. This is the subject of Theorem~\ref{orelocalisationofdg} below.

\begin{Example} \label{polynomialringisgraded} Let $K$ be a field.
Consider the polynomial algebra $K[X]$ in one variable and $X$ in degree $-1$. Then
$d(X):=1$ extends to a dg-algebra structure on $K[X]$. We get $d(X^2)=X-X=0$,
$d(X^3)=d(X^2\cdot X)=X^2=d(X\cdot X^2)$,  $d(X^4)=0=d(X^2\cdot X^2)=d(X^3\cdot X)=d(X\cdot X^3)$,
$d(X^{2n})=0$ and $d(X^{2n+1})=X^{2n}$ for all $n\in\N$.
Then $K[X]$ is $\Z$-graded, integral, and hence all non zero elements regular.
Its field of fractions is $K(X)$, the field of rational functions,
and the grading on  $K[X]$ does not extend to a grading on  $K(X)$.
\end{Example}


\begin{Prop}\label{Orelocalisationofdgrings}
Let $(R,d)$ be a differential graded ring and let $S$ be a multiplicative
subset of $\ker(d)$ {\red satisfying the Ore condition in $R$}.
Then $\textup{ass}(R):=\{r\in R\;|\;\exists s\in S\;:\;rs=0\}$ is a two-sided
dg-ideal of $(R,d)$.
\end{Prop}

Proof.
By McConnell and Robson~\cite[Section 2.1.9]{McconnellRobson} we get that $\textup{ass}(S)$ is a two-sided ideal.
We need to show that it is a dg-ideal. For this, let $r\in\textup{ass}(S)$ and let $rs=0$ for some $s\in S$.
Then for all homogeneous $r\in R$ we get
$$0=d(0)=d(rs)=d(r)s+(-1)^{|r|}r\cdot d(s)=d(r)s$$
since $S\subseteq\ker(d)$. Hence
$\textup{ass}(S)$ is a two-sided dg ideal. \dickebox
 
\begin{Rem}
The fact that we need that $S\subseteq \ker(d)$ in order to show that $ass(S)$ is a dg-ideal
is disturbing. For the sequel we shall develop the theory for homogeneous regular elements $S$. 
\end{Rem}

\begin{Rem}
Let $S$ be a multiplicative system of homogeneous regular elements of $(R,d)$. Then
$(R,d)$ is actually either unbounded or $S$ is concentrated in degree $0$. Indeed,
the $k$-th power of an element $x$ in degree $2n$ is in degree $2nk$. If $\ol x\neq 0$, then
$\ol x$ regular implies that $x^k\neq 0$ and therefore also the degree $2nk$ of $R$ is non zero.
\end{Rem}

\begin{Rem}\label{localisationatregulars}
Recall from \cite[end of 2.1.16]{McconnellRobson} that in an Ore localisation
$R_S$ we have $(a,s)=(b,t)$ for homogeneous elements $a,b,s,t$ if and only if there are
$c_1\in S$ and $a_2\in R$ such that $c_1b=a_2a\in R$ and $c_1t=a_2s\in S$.
Hence, if $S$ only contains regular elements, then the natural ring homomorphism
\begin{eqnarray*}
R&\stackrel{\lambda}\lra&R_S\\
r&\mapsto&r\cdot 1^{-1}
\end{eqnarray*}
is injective. Indeed, $\lambda(r_1)=\lambda(r_2)$ if and only if
there are $c_1\in S,c_2\in R$ such that $1\cdot c_1=1\cdot c_2\in S$ and $r_1c_1=r_2c_2$.
Then, $S$ being a set of regular elements, this implies that $r_1=r_2$.
\end{Rem}

\subsection{Localisation of dg-rings at homogeneous elements}

We prove that Ore localisation of dg-rings at homogeneous regular elements
is again a dg-ring. This is the main result of the paper.

\begin{Theorem}\label{orelocalisationofdg}
Let $(R,d)$ be a dg-ring, and let $S$ be a multiplicative set of homogeneous elements. Assume that either 
$S$ consists of regular
elements, or else $S\subseteq \ker(d)$ {\red is an  Ore set} and the image of $S$ in $R/ass(S)$ consists of regular elements of $R/ass(S)$.
Then
$$d(b,s):=(-1)^{|s|+1}(d(s),s)\cdot (b,s)+(-1)^{|s|}(d(b),s)$$
defines a differential graded structure  on $R_S$,
and the natural homomorphism
is a dg ring homomorphism $\lambda:(R,d)\lra (R_S,d_S)$ such that $\lambda(S)\in R_S^\times$
and such that for any $q\in R_S$ there is $s\in S$ with $q\cdot\lambda(s)\in\im(\lambda)$.
\end{Theorem}

Proof.
If $ass(S)\neq 0$, then by the hypothesis, $S\subseteq\ker(d)$, and using 
Proposition~\ref{Orelocalisationofdgrings},
we may replace $R$ by $R/ass(S)$, and then assume that $ass(S)=0$.

Let $t\in S$ be homogeneous. Then by the Ore condition there are
$g_1\in S$ and $a_1\in R$ such that $g_1a=a_1t$ and
for $s,t\in S$ there is $s_1,t_1\in S$ with $t_1s=s_1t$. Define then
$$(a,s)\cdot (b,t)=(a_1b,g_1s)\textup{ and }(a,s)+(b,t)=(t_1a+s_1b,t_1s).$$
As we have seen in the proof of Proposition~\ref{Orelocalisationofdgrings}
we get that $d(1)=0$. Assume that we may find a dg-ring as in the statement of the theorem
in which $s$ is invertible. Then
$$0=d(s\cdot s^{-1})=d(s)\cdot s^{-1}+(-1)^{|s|}s\cdot d(s^{-1})$$
and hence we need to define $$d(s^{-1})=(-1)^{|s|+1}s^{-1}\cdot d(s)\cdot s^{-1}.$$
Therefore, 
we can determine a
general formula for the differential.
\begin{align*}
d(b,s)=&d((1,s)\cdot (b,1))\\
=&d(1,s)\cdot (b,1)+(-1)^{|s|}(1,s)\cdot d(b,1)\\
=&(-1)^{|s|+1}(1,s)\cdot d(s,1)\cdot(1,s)\cdot (b,1)+(-1)^{|s|}(1,s)\cdot (d(b),1)\\
=&(-1)^{|s|+1}(1,s)\cdot(d(s),1)\cdot(1,s)\cdot (b,1)+(-1)^{|s|}(d(b),s)\\
=&(-1)^{|s|+1}(d(s),s)\cdot(b,s)+(-1)^{|s|} (d(b),s)
\end{align*}
Let $b\in R$ and $s\in S$.
Consider now an element $(sb,s)=(b,1)$. Then, since $d$ should extend the differential on
$R$, we get $d(b,1)=(d(b),1)$.
Then
\begin{align*}
d(sb,s)
=&(-1)^{|s|+1}(d(s),s)\cdot(sb,s)+(-1)^{|s|}(d(sb),s))\\
=&(-1)^{|s|+1}(d(s),s)\cdot(b,1)+(-1)^{|s|}(d(s)b+(-1)^{|s|}sd(b),s),s)\\
=&(-1)^{|s|+1}(d(s),s)\cdot(b,1)+(-1)^{|s|}(d(s)b,s)+(d(b),1)\\
=&(-1)^{|s|+1}(d(s)b,s)+(-1)^{|s|}(d(s)b,s)+(d(b),1)\\
=&(d(b),1)
\end{align*}
We assume that $(a,s)=(b,t)$ for homogeneous elements $a,b,s,t$. Then there are
$c_1\in S$ and $a_2\in R$ such that $c_1b=a_2a\in R$ and $c_1t=a_2s\in S$.
Then for any $t\in R$ with $ts\in S$ we get
\begin{align*}
d(ta,ts)-d(a,s)=&(-1)^{|ts|+1}(d(ts),ts)\cdot(ta,ts)+(-1)^{|ts|} (d(ta),ts)-d(a,s)\\
=&(-1)^{|ts|+1}(d(t)s+(-1)^{|t|}td(s),ts)\cdot(ta,ts)+\\
&(-1)^{|ts|} (d(t)a+(-1)^{|t|}td(a),ts)-d(a,s)\\
=&(-1)^{|ts|+1}(d(t)s,ts)\cdot(ta,ts)+(-1)^{|s|+1}(d(s),s)\cdot(ta,ts)+\\
&(-1)^{|ts|} (d(t)a,ts)+(-1)^{|s|}(d(a),s)-d(a,s)\\
=&(-1)^{|ts|+1}(d(t)s,ts)\cdot(ta,ts)+(-1)^{|ts|} (d(t)a,ts)+\\
&(-1)^{|s|+1}(d(s),s)\cdot(a,s)+(-1)^{|s|}(d(a),s)-d(a,s)\\
=&(-1)^{|ts|+1}(d(t)s,ts)\cdot(ta,ts)+(-1)^{|ts|} (d(t)a,ts)\\
=&(-1)^{|ts|}((d(t)a,ts)-(d(t)s,ts)\cdot(ta,ts))\\
=&(-1)^{|ts|}((d(t)a,ts)-((d(t),ts)\cdot(s,1)\cdot(ta,ts)))\\
=&(-1)^{|ts|}((d(t)a,ts)-((d(t),ts)\cdot(s,1)\cdot(a,s)))\\
=&(-1)^{|ts|}((d(t)a,ts)-((d(t),ts)\cdot(sa,s1)))\\
=&(-1)^{|ts|}((d(t)a,ts)-((d(t),ts)\cdot(a,1)))\\
=&(-1)^{|ts|}((d(t)a,ts)-(d(t)a,ts))\\
=&0
\end{align*}
Hence
\begin{align*}
d(a,s)=&d(a_2a,a_2s)\\
=&d(c_1b,c_1t)\\
=&d(b,t)
\end{align*}
This shows that the above definition is well-defined.

We need to verify the Leibniz formula.
We need to verify that $$d((a,s)\cdot (b,t))=d(a,s)\cdot (b,t)+(-1)^{|a|-|s|}(a,s)\cdot d(b,t)$$
for homogeneous elements $a,b\in R$ and $s,t\in S$.


We come to the general case. Let $a_1\in R$ and $g\in S$ such that $a_1t=ga$.
Let us compute the left hand term
\begin{align*}
d((a,s)\cdot (b,t))=&d((a_1b,gs))\\
=&(-1)^{|gs|+1}(d(gs),gs)(a_1b,gs)+(-1)^{|gs|}(d(a_1b),gs)\\
=&(-1)^{|gs|+1}\left[(d(g)s,gs)(a_1b,gs)+(-1)^{|g|}(gd(s),gs)(a_1b,gs)\right]\\
&+(-1)^{|gs|}\left[(d(a_1b,gs)+(-1)^{a_1|}(a_1d(b),gs)\right]\\
=&(-1)^{|gs|+1}(d(g)s,gs)(a_1b,gs)+(-1)^{|s|+1}(gd(s),gs)(a_1b,gs)\\
&+(-1)^{|gs|}(d(a_1)b,gs)+(-1)^{|gs|+|a_1|}(a_1d(b),gs)\\
=&(-1)^{|gs|+1}(d(g)s,gs)(a_1b,gs)+(-1)^{|s|+1}(d(s),s)(a_1b,gs)\\
&+(-1)^{|gs|}(d(a_1)b,gs)+(-1)^{|gs|+|a_1|}(a_1d(b),gs)
\end{align*}
The right hand term reads as
\begin{align*}
d((a,s))(b,t)+(-1)^{|(a,s)|}(a,s)d(b,t)=&(-1)^{|s|+1}(d(s),s)(a,s)(b,t)+(-1)^{|s|}(d(a),s)(b,t)\\
&+(-1)^{|(a,s)|}(a,s)\left[(-1)^{|t|+1}(d(t),t)(b,t)+(-1)^{|t|}(d(b),t)\right]\\
=&(-1)^{|s|+1}(d(s),s)(a,s)(b,t)+(-1)^{|s|}(d(a),s)(b,t)\\
&+(-1)^{|(a,s|)+|t|+1}(a,s)(d(t),t)(b,t)+(-1)^{|(a,s)|+|t|}(a,s)(d(b),t)\\
=&(-1)^{|s|+1}(d(s),s)(a_1b,gs)+(-1)^{|s|}(d(a),s)(b,t)\\
&+(-1)^{|(a,s|)+|t|+1}(a_1d(t),gs)(b,t)+(-1)^{|(a,s)|+|t|}(a,s)(d(b),t)
\end{align*}
The second term of the left hand side equals the first term of the right hand side.
We hence need to verify
\begin{align*}
(-1)^{|gs|+1}(d(g)s,gs)(a_1b,gs)&
+(-1)^{|gs|}(d(a_1)b,gs)+(-1)^{|gs|+|a_1|}(a_1d(b),gs)\\
\stackrel{!}=&
(-1)^{|s|}(d(a),s)(b,t)
+(-1)^{|(a,s|)+|t|+1}(a_1d(t),gs)(b,t)\\&+(-1)^{|(a,s)|+|t|}(a,s)(d(b),t)
\end{align*}
This is equivalent to
\begin{align*}
(-1)^{|gs|+1}(d(g)s,gs)(a_1b,gs)&
+(-1)^{|gs|}(d(a_1)b,gs)+(-1)^{|gs|+|a_1|}(a_1d(b),gs)\\
\stackrel{!}=&
(-1)^{|s|}(d(a),s)(b,t)
+(-1)^{|(a,s|)+|t|+1}(a_1d(t),gs)(b,t)\\
&+(-1)^{|(a,s)|+|t|}(a_1d(b),gs).
\end{align*}
Note that we did not assume that the elements of $S$ are regular. However, by
\cite[8.1.1 Lemma]{NastasescuVanOystaen} we see that in the equation $ga=a_1t$,
we may assume that also $a_1$ is homogeneous. Further we get $|ga_1|=|a_1t|$
and therefore $(-1)^{|a_1|+|g|}=(-1)^{|a|+|t|}$ and hence the last terms
of the left hand side and the right hand side coincide.
Therefore the equation we need to verify is equivalent to
\begin{align*}
(-1)^{|gs|+1}(d(g)s,gs)(a_1b,gs)&
+(-1)^{|gs|}(d(a_1)b,gs)\\
\stackrel{!}=&
(-1)^{|s|}(d(a),s)(b,t)
+(-1)^{|(a,s|)+|t|+1}(a_1d(t),gs)(b,t).
\end{align*}
Further, both sides are right multiples of $(b,1)$, hence
\begin{align*}
\left((-1)^{|gs|+1}(d(g)s,gs)(a_1,gs)\right.&
\left.+(-1)^{|gs|}(d(a_1),gs)\right)\cdot (b,1)\\
\stackrel{!}=&
\left((-1)^{|s|}(d(a),s)(1,t)
+(-1)^{|(a,s|)+|t|+1}(a_1d(t),gs)(1,t)\right)\cdot (b,1).
\end{align*}
Since we need to verify the equation for all $b$, we actually need to prove
\begin{align*}
(-1)^{|gs|+1}(d(g)s,gs)(a_1,gs)&
+(-1)^{|gs|}(d(a_1),gs)\\
\stackrel{!}=&
(-1)^{|s|}(d(a),s)(1,t)
+(-1)^{|(a,s|)+|t|+1}(a_1d(t),gs)(1,t).
\end{align*}
We may multiply with $(t,1)$ from the right to obtain
\begin{align*}
(-1)^{|gs|+1}(d(g)s,gs)(a_1t,gs)&
+(-1)^{|gs|}(d(a_1)t,gs)\\
\stackrel{!}=&
(-1)^{|s|}(d(a),s)
+(-1)^{|(a,s|)+|t|+1}(a_1d(t),gs).
\end{align*}
Since $ga=a_1t$ we need to show
\begin{align*}
(-1)^{|gs|+1}(d(g)s,gs)(ga,gs)&
+(-1)^{|gs|}(d(a_1)t,gs)\\
\stackrel{!}=&
(-1)^{|s|}(d(a),s)
+(-1)^{|(a,s|)+|t|+1}(a_1d(t),gs).
\end{align*}
which is equivalent to
\begin{align*}
(-1)^{|gs|+1}(d(g)s,gs)(a,s)&
+(-1)^{|gs|}(d(a_1)t,gs)\\
\stackrel{!}=&
(-1)^{|s|}(d(a),s)
+(-1)^{|(a,s|)+|t|+1}(a_1d(t),gs).
\end{align*}
This is equivalent to
\begin{align*}
(-1)^{|gs|+1}(d(g)s,gs)(a,s)&-(-1)^{|s|}(d(a),s)\\
\stackrel{!}=&
(-1)^{|(a,s|)+|t|+1}(a_1d(t),gs)-(-1)^{|gs|}(d(a_1)t,gs).
\end{align*}
However, $ga=a_1t$ implies
$$d(g)a+(-1)^{|g|}gd(a)=d(ga)=d(a_1t)=d(a_1)t+(-1)^{|a_1|}a_1d(t)$$
Hence, we need to show
\begin{align*}
(-1)^{|gs|+1}(d(g)s,gs)(a,s)&-(-1)^{|s|}(d(a),s)\\
\stackrel{!}=&
(-1)^{|(a,s|)+|t|+1}(a_1d(t),gs)-(-1)^{|gs|}(d(g)a,gs)-(-1)^{|s|}(d(a),s)\\
&+(-1)^{|a_1|+|gs|}(a_1d(t),gs).
\end{align*}
This first and the last term of the right hand side are opposite to each other, so that we need
to show
$$
(-1)^{|gs|+1}(d(g)s,gs)(a,s)-(-1)^{|s|}(d(a),s)
\stackrel{!}=
-(-1)^{|gs|}(d(g)a,gs)-(-1)^{|s|}(d(a),s).
$$
The last terms of the left and the right hand side coincide, which shows that we need to verify
$$
(-1)^{|gs|+1}(d(g)s,gs)(a,s)
\stackrel{!}=
-(-1)^{|gs|}(d(g)a,gs).
$$
The signs cancel and hence we need to show
$$
(d(g)s,gs)(a,s)\stackrel{!}=(d(g)a,gs).
$$
Since this has to hold for any $a$, we need to show
$$
(d(g)s,gs)(1,s)\stackrel{!}=(d(g)1,gs).
$$
Since we are dealing with Ore sets $S$ all this is equivalent to the ring theoretic
equation
$$s^{-1}g^{-1}d(g)ss^{-1}=s^{-1}g^{-1}d(g)$$
in the quotient ring of $R$. But this is trivially true. \dickebox

\subsection{Application: Goldie's theorem}

We have seen that Goldie's theorem~\ref{goldiestheoremclassical} makes use of
Ore localisations. We want to find a version of Goldie's theorem for differential graded rings
and take a few lines about some considerations in this direction.
We consider the case of a differential graded ring $(R,d)$.

In the set $S$ of regular elements we
consider the subset $S_{dg}$ of homogeneous elements.
By Theorem~\ref{orelocalisationofdg}  we may localise at the subset
$S_{dg}$ to get a differential graded algebra for the localisation.
The universal property of Ore localisation gives a
unique ring homomorphism
$$R_{S_{dg}}\lra Q=R_S.$$

Note however that this homomorphism is just a ring homomorphism and hence the image is not
necessarily a dg-ring.

By Remark~\ref{localisationatregulars} we get that the natural morphisms give injective maps
$$
\xymatrix{R\ar@{^{(}->}[r]& R_{S_{dg}}\ar@{^{(}->}[r]& Q=R_S.}
$$
The first map is a homomorphism of dg-rings, whereas the second is a homomorphism of rings.

\begin{Rem} Recall Example~\ref{polynomialringisgraded}.
Note that for a $\Z$-graded ring $R$ the localisation $R_S$ of $R$ at a multiplicative
subset $S$ of regular elements is $\Z$-graded (with $R\lra R_S$ a homomorphism of graded rings)
only when $S$ is homogeneous. Hence Theorem~\ref{orelocalisationofdg} is optimal.
\end{Rem}

We have an easy consequence.

\begin{Cor}
Let $(R,d)$ be a differential graded ring. If $R$ is
a $\Z$-graded $\Z$-prime right $\Z$-graded Goldie ring then the localisation of $R$
at homogeneous regular elements is a differential graded dg-simple,
dg-Artinian ring.
\end{Cor}

Proof. Apply Theorem~\ref{GoodearStafford} to the case $G=\Z$
and use that
by Theorem~\ref{orelocalisationofdg}
the localisation at regular homogeneous elements is differential graded.

\subsection{Comparing homology localisation with Ore localisation}

Let $(R,d)$ be a differential graded ring.

Let $S$ be a multiplicative system of homogeneous elements of even degree in $\ker(d)$ such
that the image of $S$ in $R/ass(S)$
contains only regular elements, then, following Proposition~\ref{Orelocalisationofdgrings}
we may form the Ore localisation $R_S$ at $S$. However, since $S\subseteq \ker(d)$, we may form
the image $\ol S$ in $H(R,d)$. Since $S$ is a multiplicative system, $\ol S$ is a multiplicative
system in $H(R,d)$. Suppose for the moment that $(R,d)$ is a differential graded $k$-algebra for a field $k$.
Then, by \cite[Theorem 3.10]{BraunChuangLazarev} we see that
$$R\ast_{k(\ol S)}^{\mathbb L}k(\ol S,{\ol S}^{-1})$$
is the universal ring inverting all elements of $\ol S$. Here,
$\ast_{k(\ol S)}^{\mathbb L}k(\ol S,{\ol S}^{-1})$
denotes the derived coproduct in the category of dg-rings. It can be computed
by replacing the left hand argument with a cofibrant dg-algebra under $k(\ol S)$.
More precisely (cf \cite[Definition 2.4]{BraunChuangLazarev}), consider the under
category $A\dar dgAlg$ formed by objects being
dg-algebra homomorphisms $A\lra C$ and morphisms being commutative triangles.
For any dg-algebra homomorphism $A\lra B$ we obtain the restriction functor
$$B\dar dgAlg\lra A\dar dgAlg$$
and see that this has a left adjoint denoted by $B\ast_A-$. The derived functor, replacing $A$ and $B$ by
cofibrant replacements, is then denoted by $B\ast_A^{\mathbb L}-$.
Hence, there is a unique homomorphism
of dg $k$-algebras $$R\ast^{\mathbb L}_{k(\ol S)}k(\ol S,{\ol S}^{-1})\stackrel{\lambda}{\lra} R_S$$
such that the diagram
$$
\xymatrix{
R\ast^{\mathbb L}_{k(\ol S)}k(\ol S,{\ol S}^{-1})\ar[rr]^-{\lambda}&& R_S\\
&R\ar[lu]\ar[ru]
}
$$
is commutative. However, since $R\ast_{k(\ol S)}^{\mathbb L}k(\ol S,{\ol S}^{-1})$
does not necessarily invert all elements of $S$, but only those in $\ol S$, we do not necessarily get that
$\lambda$ is invertible. An example of this kind occurs if $(R,d)$ is acyclic.

\medskip

\begin{Prop}\label{BCLisquasiiso}
Suppose that $S\subseteq\ker(d)$ is a multiplicative {\red Ore} set of homogeneous elements, 
and the image of $S$ in $R/ass(S)$ consists of regular elements in $R/ass(R)$. Then
$H(\lambda)$ is an isomorphism and hence $R\ast_{k(\ol S)}^{\mathbb L}k(\ol S,{\ol S}^{-1})$
is quasi-isomorphic to $R_S$. In particular, if $\ol S$ is the image of $S$ in $H(R)$, then
$$H(R_S)\simeq H(R)_{\ol S}.$$
\end{Prop}

Proof.
We consider $H(R_S)$. By \cite[Proposition 5.14]{BraunChuangLazarev}
we get $H(R_S)=H(R)_{\ol S}$. This then shows that $H(\lambda)$ is an isomorphism and hence
$\lambda$ is a quasi-isomorphism since both
$H(R\ast^{\mathbb L}_{k(\ol S)}k(\ol S,{\ol S}^{-1}))$
and $H(R)_{\ol S}$ have the same universal property being initial with respect to inverting $\ol S$.
\dickebox

\begin{Rem}
Note that the authors of \cite{BraunChuangLazarev} mention that their construction 
of $R\ast_{k(\ol S)}^{\mathbb L}k(\ol S,{\ol S}^{-1})$ is hard to perform explicitly. 
Our Proposition~\ref{BCLisquasiiso} provides such a construction up to quasi-isomorphism
and our localisation is a lot more general. 
\end{Rem}

\section{Semiprime differential graded rings}

\label{dgprimesect}

We want to prove a dg-version of the results of Section~\ref{recall}.

\subsection{dg-prime, dg-semiprime}

\begin{Def}
Let $(R,d)$ be a differential graded ring. We consider nilpotent differential graded ideals $(I,d)$
and define {\em $\dgnil(R,d)$ to be the sum of all nilpotent differential graded
nilpotent two-sided ideals}.
\end{Def}

\begin{Def}
A dg-ring $(R,d)$ is called {\em dg-Noetherian} if any ascending chain of dg-ideals
of $(R,d)$ is finite. A dg-ring $(R,d)$ is  {\em dg-Artinian} if any descending chain of dg-ideals
of $(R,d)$ is finite.
\end{Def}

\begin{Lemma}\label{dgnilisdg-ideal}
If $(R,d)$ is a dg-ring, then $\dgnil(R,d)$ is a dg-ideal of $(R,d)$, which we call
the {\em differential graded nil radical.} If $(R,d)$ is dg-Noetherian, then $\dgnil(R,d)$ is nilpotent.
\end{Lemma}

Proof.
Since the sum of ideals is an ideal, and since the differential is additive, and elements of a
sum of ideals is a finite sum of elements of the constituants, the sum of dg-ideals is a dg-ideal.

Since the sum of two nilpotent ideals is a nilpotent ideal, since the sum of two differential
graded ideals is a differential graded ideal, at least if $(R,d)$ is dg-noetherian,
using the ascending chain condition, the sum of all nilpotent differential
graded nilpotent two-sided ideals is actually a finite sum of
nilpotent differential graded nilpotent two-sided ideals.
Hence if $(R,d)$ is dg-noetherian, then $\dgnil(R,d)$ is a nilpotent differential
graded ideal. \dickebox

\begin{Lemma}\label{dgnillemma}
Let $(R,d)$ be a dg-Noetherian differential graded ring. Then $\dgnil(R/\dgnil(R,d))=0$.
\end{Lemma}

Proof. 
Since Lemma~\ref{dgnilisdg-ideal} shows that $\dgnil(R,d)$ is a dg-ideal, and since
any quotient $R/I$ of a differential graded algebra by a differential
graded ideal $(I,d)$ is again a differential
graded ring, $d$ induces a differential $\ol d$ on $R/\dgnil(R)$.
Also, since preimages of dg-ideals under dg-homomorphisms are dg-ideals,
$(R/\dgnil(R,d),\ol d)$ is dg-Noetherian. Hence, if $(\ol I,\ol d)$ is a nilpotent
differential graded ideal of $(R/\dgnil(R,d),\ol d)$, its preimage $I$ is a differential graded
ideal $(I,d)$. Since there is an integer $n$ such that $\ol I^n=0$, and hence
$(I,d)^n\subseteq\dgnil(R,d)$. Now, $\dgnil(R,d)^m=0$ for some $m$, since $(R,d)$ is dg-Noetherian,
and therefore $I^{n+m}=0$. This then shows $I\subseteq\dgnil(R,d)$ and hence $\ol I=0$.
\dickebox

\begin{Def}
Let $(R,d)$ be a differential graded ring.
\begin{itemize}
\item A two-sided differential graded ideal $(P,d)$ is called {\em dg-prime} if
whenever $(S,d)$ and $(T,d)$ are
two-sided differential graded ideals with $ST\subseteq P$, then $S\subseteq P$ or $T\subseteq P$.
\item $(R,d)$ is called {\em dg-semiprime} if $\dgnil(R)=0$.
\item $(R,d)$ is called {\em dg-prime} if for all non zero two-sided dg-ideals
$(I,d)$ and $(J,d)$ we get $IJ\neq 0$.
\end{itemize}
\end{Def}

Again, if $R$ is concentrated in degree $0$ (and $d=0$), then the concept of dg-(semi-)prime
coincides with the concept of (semi-)prime.

\begin{Lemma}
dg-Noetherian dg-prime rings $(R,d)$ are dg-semiprime.
\end{Lemma}

Proof. Indeed, $\dgnil(R,d)$ is nilpotent, say
$\dgnil(R,d)^k=0$. Then $$\dgnil(R,d)\cdot \dgnil(R,d)^{k-1}=0.$$ Since $(R,d)$ is assumed
to be semiprime, we get $\dgnil(R,d)^{k-1}=0$ (or $\dgnil(R,d)=0$ which implies the former)
and by induction on $k$ we get $\dgnil(R,d)=0$, whence $(R,d)$ is dg-semiprime.
\dickebox

\medskip

Let $R$ be a Noetherian ring.
Recall that the nil radical $\Nil(R)$ of a ring $R$ is the sum of all nilpotent
two-sided ideals of $R$.
It is, by definition, the largest nilpotent ideal of $R$.
Further, it is a classical result that $\Nil(R)$ is the intersection of all prime ideals of $R$.

\begin{Example}\label{dgPraduneaqualdgnil}
Let $K$ be a field and let $A=K[X]/X^2$. Then $A$ is graded when we declare $X$ to be in degree $-1$.
Further, $d(X)=1$ and $d(1)=0$ gives a structure of differential graded algebra on $A$.
The only ideals are $0$, $XK[X]/X^2$, and $A$. The ideal $XK[X]/X^2$ is not differential graded
and hence the intersection of dg-prime ideals of $A$ {\red is $0$, as well as} 
$\dgnil(A,d)=0$. Note that the classical prime radical is $XA$, hence larger. 
{\red Recall the classical result that for artininan algebras} $\Nil(R)$ is the intersection of all
prime ideals.
\end{Example}

\begin{Def}
Let $(R,d)$ be a differential graded ring. Then the {\em dg-prime radical}
$\textup{Prad}(R,d)$ is the intersection of all dg-prime ideals of $(R,d)$.
\end{Def}

\begin{Lemma} \label{PradofRmodPrad}
Let $(R,d)$ be a differential graded ring. Then $\Prad((R/\Prad(R,d),\ol d))=0$.
\end{Lemma}

Proof. Let $\ol I:=\Prad((R/\Prad(R,d),\ol d))$ and denote by $I$ the preimage of $\ol I$ in $R$.
Then clearly $\Prad(R,d)\subseteq I$. We show that the inclusion in the opposite sense also holds.

In order to do so, we need to show that $I$ is contained in every dg-prime ideal $P$ of $R$.
Let $Q$ be such a dg-prime ideal. Then the image $\ol Q$ of $Q$ in $R/\Prad(R,d)$
is again a dg-prime ideal of $R/\Prad(R,d)$. Hence $\ol I\subseteq\ol Q$, and
therefore $I\subseteq Q$.
This shows the lemma. \dickebox

\begin{Lemma}\label{dgnilisinprad}
Let $(R,d)$ be a  differential graded ring.
Then $\dgnil(R,d)\subseteq \Prad(R,d)$.
\end{Lemma}

Proof. Indeed, let $L$ be a nilpotent differential graded ideal of $R$. Then
$L^k=0$ for some $k$. Let $Q$ be a dg-prime ideal of $(R,d)$. Then
$L\cdot L^{k-1}=L^k=0\subseteq Q$ and hence, since $Q$ is dg-prime,
$L^{k-1}\subseteq Q$. By induction on $k$ we get that $L\subseteq Q$. Hence $L\subseteq \Prad(R,d)$
since $Q$ is arbitrary dg-prime. Since $\dgnil(R,d)$ is the sum of all nilpotent
differential graded ideals, we also get $\dgnil(R,d)\subseteq\Prad(R,d)$.
This shows the lemma. \dickebox

\begin{Lemma}
Let $(R,d)$ be a dg-Noetherian differential graded ring.
Then $(R/\dgnil(R,d),\ol d)$ is dg-semiprime.
\end{Lemma}

Proof.
This is an immediate consequence of Lemma~\ref{dgnillemma}. \dickebox

\begin{Lemma}\label{modPradissemiprime}
Let $(R,d)$ be a dg-Noetherian differential graded ring.
Then $(R/\Prad(R,d),\ol d)$ is dg-semiprime.
\end{Lemma}

Proof.
By Lemma~\ref{dgnilisinprad} we get a surjective homomorphism of differential graded rings
$$(R/\dgnil(R,d),\ol d_1)\lra (R/\Prad(R,d),\ol d_2)$$
given by the natural inclusion from  Lemma~\ref{dgnilisinprad}.

Let $\ol I:=\dgnil(R/\Prad(R,d),\ol d_2)$, and let $I$ be the preimage of $\ol I$ in
$R/\dgnil(R,d)$. Since $R$ is dg-Noetherian, also $R/\Prad(R,d)$ is dg-Noetherian, and hence
$\ol I$ is nilpotent. Therefore there is $k\in\N$ such that $I^k\subseteq \Prad(R,d)/\dgnil(R,d)$.
But this implies that $I\subseteq \Prad(R,d)/\dgnil(R,d)$ by the defining property of a
dg-prime ideal, and a standard induction on $k$.
Hence $\ol I=0$ and we proved the lemma. \dickebox



\begin{Def}
Let $(R,d)$ be a differential graded ring. We say that $(R,d)$ is {\em strongly
dg-semiprime} if $\Prad(R,d)=0$.
\end{Def}


{\red 
Just like in the classical case we get analogously the following

\begin{Prop}\label{dgprimeradicalforartinian}
Let $(A,d)$ be a differential graded algebra and suppose that $(A,d)$ is dg-artinian on 
twosided differential graded ideals.
Then $\Prad((A,d)=\dgnil(A,d)={\textup {dgrad}}_2(A,d)$. 
\end{Prop}

Proof. 
In a first step we show that $(A,d)$ only contains a finite number of maximal twosided 
differential graded ideals. Indeed, if $I_1,I_2,\dots$ is a sequence of maximal 
twosided differential graded ideals, then 
$$I_1\supsetneq I_1\cap I_2\supsetneq\dots$$
is a strictly descending sequence of twosided differential graded ideals of $(A,d)$. 
Hence this has to stop, by the dg-artinian property on twosided dg-ideals. This shows the first step. 

\medskip

In a second step we show that ${\textup {dgrad}}_2(A,d)$ is nilpotent. Indeed, this is a direct consequence of the dg-Nakayama Lemma~\cite[Lemma 4.28]{dgorders} and the dg-Artinianity on
twosided dg-ideals. 

\medskip

The third step shows that any dg-prime ideal is twosided maximal. Indeed, if $\wp $ is
a dg-prime ideal, then since  ${\textup {dgrad}}_2(A,d)$ is nilpotent, also 
${\textup {dgrad}}_2(A,d)\subseteq\wp$. But if $I_1,\dots,I_n$ are the (finite number!) 
maximal twosided differential graded ideals, then 
$$I_1\cdot I_2\cdot\cdot\dots\cdot I_n\subseteq {\textup {dgrad}}_2(A,d)\subseteq\wp$$
and since $\wp$ is dg-prime, there is $j$ such that $I_J\subseteq\wp.$ 
Since $I_j$ is maximal, we have equality.

We proved the proposition.
\dickebox
}

\section{dg-essential dg-submodules}

\label{dgessential}

Recall from \cite[Section 2.3.4]{NastasescuVanOystaen} that for a $G$-graded ring $R$
and a $G$-graded $R$-module $M$ we say that a submodule $N$ of $M$ is
gr-essential if for any $G$-graded submodule $X$ of $M$ we have $X\cap N\neq 0$.

\begin{Def}
A non zero submodule  $M$ of a differential graded
module $(X,\delta)$ over a differential graded ring $(R,d)$ is called {\em dg-essential} if
for any differential graded submodule $(N,\delta)$ of $(X,\delta)$ one has $M\cap N\neq 0$.
A {\em dg-essential ideal} is a dg-essential submodule of the regular module $(R,d)$.
A {\em dg-essential two-sided ideal} is a dg-essential submodule of the
$(R\otimes_\Z R^{op},d\otimes d^{op})$-module $(R,d)$.
\end{Def}

Note that for a submodule $M$ of $(X,\delta)$
we did not assume that $M$ is differential graded. We shall need this
subtlety later. However, most of the time we shall assume that $M$ is differential graded.

If $X$ is ungraded and $\delta=0$, we get back the usual concept of an essential submodule.
Similarly, if the grading is non zero, but the differential is zero, we get the concept
of a gr-essential submodule (cf e.g. \cite{NastasescuVanOystaen}).

\begin{Lemma}\label{gressentialimpliesdgessential}
Further, let $(R,d)$ is a dg-ring, and let $(M,\delta)$ be a dg-$(R,d)$-module.
If $(N,\delta)$ is a dg-submodule, and if the $\Z$-graded submodule $N$ is gr-essential
in $M$, then $(N,\delta)$ is dg-essential.
\end{Lemma}

Proof. Indeed, if $(X,\delta)$ is a dg-submodule of $(M,\delta)$, then forgetting
the differential, $X$ is a $\Z$-graded submodule, and hence $X\cap N\neq 0$.
Hence $(N,\delta)$ is dg-essential in $(M,\delta)$. \dickebox

\begin{Lemma}\label{essentiallemma}
Let $(R,d)$ be a differential graded ring.
\begin{enumerate}
\item\label{essentiallemmaitem1}
If $(R,d)$ is dg-prime, then any non zero two-sided dg-ideal is dg-essential.
\item\label{essentiallemmaitem2}
The relation of being a dg-essential submodule is transitive.
\item\label{essentiallemmaitem3}
The intersection of two dg-essential submodules is dg-essential.
\item\label{essentiallemmaitem4}
If $U$ is a dg-module, and if $N$ is a dg-essential dg-module in the dg-module $M$, then
$U\oplus N$ is dg-essential in $U\oplus M$.
\item
If $(N_i,\delta_i)$ is a dg-essential submodule of $(M_i,\delta_i)$ for all $i\in\{1,\dots,n\}$, then
$(\bigoplus N_i,\bigoplus \delta_i)$ is a dg-essential submodule of  $(\bigoplus M_i,\bigoplus \delta_i)$.
\item \label{essentiallemmaitem6}
If $(N,\delta)$ is a differential graded submodule of $(M,\delta)$. Then there is a
differential graded submodule
$(X,\delta)$ of $M$ with $N\cap X=0$ and $N\oplus X$ is dg-essential in $M$.
\end{enumerate}
\end{Lemma}

Proof. \begin{enumerate}
\item
Let $(I,d)$ be a non zero differential graded two-sided ideal, and let $(X,d)$
be another non zero differential graded two-sided ideal. Then $0\neq IX\subseteq I\cap X$
since $(R,d)$ is dg-prime.
\item
If $I\leq J\leq K$ is a sequence of differential graded submodules of the dg-submodule $(M,\delta)$,
if $I$ is dg-essential in $J$, and $J$ is dg-essential in $K$,
and $(X,\delta)$ is another dg-submodule in $K$, then $X\cap J$ is a dg-submodule of $J$.
Since $J$ is dg-essential in $K$, we get $X\cap J\neq 0$.
Since $I$ is dg-essential in $J$, we get $(X\cap J)\cap I\neq 0$. Since $I\leq J$, we get $X\cap J\cap I=X\cap I$.
Hence $I$ is dg-essential in $K$.
\item
Suppose that $N_1$ is dg-essential in $(M,\delta)$ and $N_2$ is dg-essential in $(M,\delta)$, and let
 $X$ be another dg-submodule in $M$.
Since $N_1$ is
dg-essential in $M$, we get $X\cap N_1\neq 0$. Since $N_2$ is dg-essential in $M$, we get
$X\cap (N_1\cap N_2)=(X\cap N_1)\cap N_2\neq 0$. Hence $N_1\cap N_2$ is dg-essential.
\item
If $N$ is dg-essential in $M$, then $U\oplus N$ is dg-essential in $U\oplus M$. Indeed,
We denote by $\pi:U\oplus M\lra M$ the canonical projection. Let $X$ be a dg-submodule of $U\oplus M$.
Then $\pi(X)$ is a submodule of $\pi(U\oplus M)=M$. Either, $\pi(X)=0$ or, using that $N$ is
essential in $M$, we get $\pi(X)\cap \pi(U\oplus N)\neq 0$.
If $\pi(X)=0$, then $X\subseteq U$ and hence $X\cap(U\oplus N)=X\oplus 0\neq 0$.
If $\pi(X)\cap \pi(U\oplus N)\neq 0$, let $0\neq x\in \pi(X)\cap \pi(U\oplus N)$.
Then there is $u\in U$ such that $(u,x)\in X$ But then $(u,x)\in U\oplus N$, and therefore
$(U\oplus N)\cap X\neq 0$.
\item We proceed by induction on $n$. Since $N_1$ is dg-essential in $M_1$, by the previous statement
we have $N_1\oplus N_2$ is dg-essential in $M_1\oplus N_2$. Again by the previous statement
$M_1\oplus N_2$ is dg-essential in $M_1\oplus M_2$. By the second statement $N_1\oplus N_2$
is dg-essential in $M_1\oplus M_2$.
We may assume that $N_1\oplus\dots \oplus N_{n-1}$ is dg-essential in  $M_1\oplus\dots \oplus M_{n-1}$.
By the case of two factors we see that $(N_1\oplus\dots \oplus N_{n-1})\oplus N_n$ is dg-essential in
$(M_1\oplus\dots \oplus M_{n-1})\oplus M_n$.
\item
Let $\mathcal X$ be the set of dg-submodules $(Z,\delta)$ of $(M,\delta)$ such that
$(N\cap Z)=0$. Since $(0,\delta)$ is in $\mathcal X$, we get that ${\mathcal X}\neq \emptyset$.
Clearly $\mathcal X$ is partially ordered by inclusion. If $\mathcal Y$ is a totally ordered subset of
$\mathcal X$, we get
$$\widehat Y:=\bigcup_{Y\in{\mathcal Y}}Y$$
is a dg-submodule of $(M,\delta)$. Further, $\widehat Y\cap N=0,$ since else there
is $0\neq y\in \widehat Y\cap N$. Then $y\in\widehat Y$ implies that there is $Y\in{\mathcal Y}$ with
$y\in Y$. But this contradicts $Y\cap N=0$. Hence,
by Zorn's lemma there is a maximal element $(X,\delta)$ of $\mathcal X$.
By definition $N\cap X=0$.
Let $Y$ be a differential graded submodule of $M$ with $Y\cap(N\oplus X)=0$. Then $X\oplus Y$
still is a differential graded submodule satisfying
$N\cap(X\oplus Y)=0$. By maximality of $X$ we get $Y=0$. Therefore $N\oplus X$
is dg-essential in $M$.
\end{enumerate}
This proves the lemma. \dickebox

\medskip

Lemma~\ref{essentiallemma}.(\ref{essentiallemmaitem6}) suggests the following definition.

\begin{Def}
Let $(R,d)$ be a differential graded ring and let $(M,\delta)$ be a differential graded $(R,d)$-module.
For a differential graded $(R,d)$-submodule $(N,\delta)$ of $(M,\delta)$ we say that a differential graded
$(R,d)$-submodule $(L,\delta)$ is a {\em dg-complement} to $(N,\delta)$ if the following two conditions hold: 
$N\cap L=0$, and $(L,\delta)$ is maximal with respect to this property.
\end{Def}

\begin{Rem}\label{modulepluscomplementisessential}
As a consequence, if $(N,\delta)$ is a dg-submodule of $(M,\delta)$, and if $(L,\delta)$ is a
dg-complement to $(N,\delta)$ in $(M,\delta)$, then $N\oplus L$ is dg-essential. By
Lemma~\ref{essentiallemma}.(\ref{essentiallemmaitem6}) such a complement always exist.
\end{Rem}

\begin{Rem}\label{dgcomplementsneednotbegrcomplements}
Note that a dg-complement need not be a gr-complement. Moreover, a dg-essential
submodule need not be a
gr-essential submodule.
\end{Rem}

Analogous to the classical case we get

\begin{Cor}
Let $(R,d)$ be a differential graded ring. Then a differential graded $(R,d)$-module $(M,\delta)$ is
a direct sum of simple differential graded $(R,d)$-modules if and only if $0$ and $M$ are the only
dg-essential submodule of $(M,\delta)$.
\end{Cor}

Proof. If $(M,\delta)$ is a direct sum of dg-simple submodules, say $M=\bigoplus_{i\in I}M_i$,
and let $(N,\delta)$ be a dg-essential submodule of $M$. Then $N\cap M_i\neq 0$ for all $i\in I$ since $N$
is dg-essential and $(M_i,\delta)$ is a dg submodule of $M$. Since $M_i$ is dg-simple, and since $N\cap M_i\leq M_i$,
we get that $N\cap M_i=M_i$, whence $M_i\subseteq N$. This holds for all $i$, and therefore $N=M$.

Conversely, suppose that $0$ and $M$ are the only dg-essential submodules of $M$. Then,
Lemma~\ref{essentiallemma}.\ref{essentiallemmaitem6} implies that
any differential graded submodule is complemented by a differential graded submodule, and
therefore, by \cite[Lemma 4.17]{dgorders}, which is formulated for finite sums only, but which can be
generalised to arbitrary sums by the usual application of Zorn's lemma,
$(M,\delta)$ is a direct sum of simple dg-submodules as differential graded module.  \dickebox

\section{dg-uniform dimension, the dg-singular ideal}

\label{dguniformsect}

Recall that for an algebra $A$ an $A$-module $M\neq 0$ is called uniform if all non zero submodules of $M$ are
essential. A module $M$ has finite uniform dimension if $M$ does not contain an infinite direct sum of non zero submodules.
If $M$ contains a finite direct sum $\bigoplus_{i=1}^nN_i=:N$ and such that $N$ is essential, and such that
each $N_i$ is uniform, then we say that $n$ is the uniform dimension of $M$.

We can easily transpose this concept to the differential graded situation.

\begin{Def}
Let $(A,d)$ be a differential graded algebra.
\begin{itemize}
\item
A non zero differential
graded $(A,d)$-module $(M,\delta)$ is called
{\em dg-uniform}
if all non zero differential graded $(A,d)$-submodules $(N,\delta)$ of $(M,\delta)$ are dg-essential.
\item
A dg-module $(M,\delta)$ is said to have {\em finite dg-uniform dimension} if $(M,\delta)$
does not contain an infinite
direct sum of differential graded submodules.
\item
If $(M,\delta)$ contains a dg-essential
dg-submodule $(N,\delta)$ which is
the direct sum of dg-uniform submodules $N_1,\dots,N_n$, then we say that $n$ is
the {\em dg-uniform dimension of $(M,\delta)$} and write $\dgudim(M,\delta)$ for the
dg-uniform dimension, or $\dgudim_{(A,d)}(M,\delta)$ in case we need to make precise the
dg-ring which operates.
\end{itemize}
\end{Def}

As for the classical case we shall need to show that the dg-uniform dimension is well-defined.
But the proof of the classical case \cite[2.2.7, 2.2.8, 2.2.9]{McconnellRobson} carries through
verbatim. In particular,

\begin{Prop} \label{sumofuniformsubmods}
Let $(M,\delta)$ be a differential graded $(R,d)$module with finite uniform dimension. Suppose that
$\bigoplus_{i=1}^n(U_i,\delta)$ be a dg-essential submodule of $(M,\delta)$
such that each $(U_i,\delta)$ is uniform for each $i$, then any direct sum of dg-submodules
of $(M,\delta)$ has at most $n$ non zero terms, and a direct sum of non
zero dg-uniform submodules of $(M,\delta)$ is dg-essential if and only if the sum
 has $n$ terms.
\end{Prop}

Proof. Indeed, the proof of \cite[2.2.9]{McconnellRobson} carries through
verbatim. \dickebox

\medskip

Analogous to  \cite[2.2.10]{McconnellRobson} we get for the dg-situation

\begin{Lemma} \label{dgudimofsubmodules}
Let $(R,d)$ be a differential graded ring and let $(M,\delta)$, $(M_1,\delta_1)$,
$(M_2,\delta_2)$ be  differential graded $(R,d)$-modules. Then
\begin{enumerate}
\item $\dgudim(M,\delta)=1$ if and only if $(M,\delta)$ is uniform.
\item If $(N,\delta)$ is a dg-submodule of $(M,\delta)$ and $\dgudim(M,\delta)=n$, then
$\dgudim(N,\delta)\leq n$ and
$$\dgudim(N,\delta)= n\Leftrightarrow N\textup{ is dg-essential in }M$$
\item $\dgudim(M_1\oplus M_2,\delta_1\oplus\delta_2)=\dgudim(M_1,\delta_1)+\dgudim(M_1,\delta_2)$
\end{enumerate}
\end{Lemma}

Proof.
The first item is simply the definition.

For second item let $N_1,\dots,N_t$ be dg-submodules of $(N,\delta)$ such that
$N_1+\dots+N_t=N_1\oplus \cdots\oplus N_t$. Since $(N,\delta)\leq (M,\delta)$, this direct
sum $N_1\oplus \cdots\oplus N_t$ of dg-submodules of $N$ is also a direct sum of
dg-submodules of $M$. Hence $$\dgudim(N,\delta)\leq \dgudim(M,\delta).$$
If
$$\dgudim(N,\delta)= \dgudim(M,\delta),$$
then by definition $N_1\oplus \cdots\oplus N_t$ is dg-essential in $(M,\delta)$,
but then also $(N,\delta)$ is dg-essential in $(M,\delta)$ since it contains the direct sum.
If $(N,\delta)$ is dg-essential in $(M,\delta)$, and since
$N_1\oplus \cdots\oplus N_t$ is dg-essential in $(N,\delta)$, by Lemma~\ref{essentiallemma}.

The third item is trivial and follows by the definition.

This proves the Lemma. \dickebox

\begin{Rem}
Since any dg-submodule is a submodule, it is clear that $\dgudim(M,\delta)\leq \textup{udim}(M)$
for any differential graded module $(M,\delta)$. In case $M$ is concentrated in degree $0$
(and as a consequence $\delta=0$), then
$\dgudim(M,\delta)= \textup{udim}(M)$. The inequality may be strict for general $\delta\neq 0$.
We shall see an instance in Example~\ref{examplematrixring}
\end{Rem}

%
%
%

Recall from e.g. \cite[Chapter 3, \S 7]{Lam}
the classical notion of a singular module.

\begin{Def} (cf  e.g. \cite[Chapter 3, \S 7]{Lam})
Let $M$ be an $R$-right module. Then
\begin{itemize}
\item
$m\in M$ is {\em singular} if $\ann(m)=\{r\in R\;|\;mr=0\}$ is essential in $R$.
\item
The {\em singular submodule} of $M$ is the set of singular elements in $M$.
\item A module
is called singular it all elements are singular. It is
not difficult to show that this is indeed a submodule.
\item The {\em right singular ideal} is
the singular submodule of $R_R$, and the left singular ideal is the singular
submodule of $_RR$.
\item
Accordingly, for a subset $I$ of $(R,d)$ we
denote $\rann(I):=\{r\in R\;|\;Ir=0\}$ and $\lann(I):=\{r\in R\;|\;rI=0\}$.
\end{itemize}
\end{Def}

\begin{Def}
Let $(R,d)$ be a differential graded ring. Then the {\em right singular dg-ideal
$\zeta_{dg}(R,d)$} is formed by those $a\in R$ such that there is a dg-essential
differential graded $(R,d)$-right ideal $(E,d)$ with $a\cdot E=0$.
\end{Def}

As we see, we need to deal with annihilators in the dg-context. Let us give some 
elementary observations.

\begin{Lemma}\label{annihilatrosaredgideals}
Let $(R,d)$ be a differential graded ring and $I$ a subset of $R$.
\begin{enumerate}
\item
Then $\rann(I)$ is a right ideal and $\lann(I)$ is a left ideal.
\item
If $I$ is a left ideal, then $\lann(I)$ is a two-sided ideal.
\item
If $I$ is a right ideal, then $\rann(I)$ is a two-sided ideal.
\item
If $I\subseteq \ker(d)$, and $I$ is graded, then $\rann(I)$ is a dg-right ideal
and $\lann(I)$ is a dg-left ideal.
\item
If $(I,d)$ is a dg-left ideal, then $\lann(I,d)$ is a two-sided dg-ideal.
\item
If $I$ is a dg-right ideal, then $\rann(I,d)$ is a two-sided dg-ideal.
\end{enumerate}
\end{Lemma}

Proof.
\begin{enumerate}
\item
If $Ir=0$, then also $Irs=0$ for all $s\in R$ and if $Ir_1=Ir_2=0$, then
$I(r_1+r_2)\subseteq Ir_1+Ir_2=0$. Hence $\rann(I,d)$ is a right ideal.
Similarly, $\lann(I)$ is a left ideal.
\item
If $I$ is a right ideal, then $Ir=0$ implies $Isr\subseteq Ir=0$.
Hence $\rann(I)$ is a two-sided ideal.
\item
Similarly, if $I$ is a left ideal, then $\lann(I)$ is a two-sided ideal.
\item
Since $I$ is supposed to be graded, also $\rann(I)$ and $\lann(I)$ are graded.
Let $x\in \rann(I)$ be homogeneous and $z\in I$. Then
$$0=d(0)=d(zx)=d(z)x+(-1)^{|z|}zd(x)=(-1)^{|z|}zd(x)$$
and hence $d(x)\in\rann(I)$ as well. Therefore $\rann(I)$ is a dg-right ideal.
The case of $\lann(I)$ is analogous.
\item
If $(I,d)$ is a dg-right ideal, then, by the previous $\rann(I)$
is a two-sided ideal.
Further, for all homogeneous $x\in I$ and $r\in R$ we get
$$d(xr)=d(x)r+(-1)^{|x|}xd(r)$$
which implies that $d(x)r\in I$ for all $x\in I$ and $r\in R$.
Hence $\rann(I,d)$ is a two-sided dg-ideal.
\item
Similarly, if $(I,d)$ is a dg-left ideal, then $\lann(I,d)$ is a two-sided dg-ideal.
\end{enumerate}
This proves the lemma. \dickebox

\begin{Example}\label{examplematrixring}
\begin{enumerate}
\item
We recall from \cite{dgorders} the following differential graded ring. Let $R$ be any
integral domain and let
$$A=\End^\bullet_R(R\stackrel{\lambda}\lra R)$$
Then $$A=\left(\begin{array}{cc}R&R\\R&R\end{array}\right)$$
where the main diagonal is the set of degree $0$ eleemnts, the lower
diagonal is in degree $-1$ and the upper diagonal is in degree $1$.
The differential is
$$d(\left(\begin{array}{cc}0&0\\1&0\end{array}\right)):=
\left(\begin{array}{cc}\lambda&0\\0&\lambda\end{array}\right)$$
and $$d(\left(\begin{array}{cc}x&0\\0&y\end{array}\right)):=
\left(\begin{array}{cc}0&\lambda(y-x)\\0&0\end{array}\right)$$
Then $\ker(d)=R[X]/X^2$ where $X$ is in degree $1$. If $R=K$ is a field, then
$$\zeta(\ker(d))=\zeta(K[X]/X^2)=\soc(K[X]/X^2)=XK[X]/X^2).$$
Suppose from now on that $R=K$ is a field. Then the algebra $A$ is semisimple.
Right ideals correspond to rows of the matrix ring and the only non trivial
differential graded
ideal is $$\left(\begin{array}{cc}K&K\\0&0\end{array}\right)=:I$$ since it
needs to be stable by the differential and if the right lower coefficient
is non zero, then also the left lower coefficient (since it is a right ideal)
and by the differential also the upper two coefficients. Hence
$I$ is dg-essential, and it is the only dg-essential right ideal.
Note that $$\rann(\left(\begin{array}{cc}0&0\\0&1\end{array}\right))=
\rann(\left(\begin{array}{cc}0&1\\0&0\end{array}\right))=I$$
and linear combinations of these two elements are the only ones with right annihilator $I$.
Hence $$\zeta_{dg}(A,d)=\left(\begin{array}{cc}0&K\\0&K\end{array}\right).$$
This ideal is a left ideal only, and in particular is not a twosided ideal.
Further,
$$\zeta_{dg}(A)\cdot\zeta_{dg}(A)=\left(\begin{array}{cc}0&K\\0&K\end{array}\right)\cdot \left(\begin{array}{cc}0&K\\0&K\end{array}\right)=
\left(\begin{array}{cc}0&K\\0&K\end{array}\right)=\zeta_{dg}(A).$$
Hence, $\zeta_{dg}(A)$ is not nilpotent, unlike \cite[Lemma 3.4]{McconnellRobson}
in the classical case.
Since the algebra $A$ is semisimple as an algebra, $\zeta(A)=0$.
\item
Let $K$ be a field. Then $A=K[X]/X^2$ is a dg-algebra with $d(X)=1$.
Then there is no non trivial dg-ideal of $A$. Hence, $\zeta_{dg}(A)=0$. However,
$XA=\soc(A)$ and hence any ideal intersects non trivially with $\soc(A)$. Furthermore
$\rann(X)=\soc(A)$. Therefore $\zeta(A)=\soc(A)$.
\end{enumerate}
\end{Example}

\begin{Prop}\label{singularidealofkernel}
Let $(R,d)$ be a dg-ring.  Then
$$\zeta(R)\cap\ker(d)\subseteq \zeta_{dg}(R,d)\cap\ker(d)$$
and
$$\zeta(\ker(d))\subseteq\zeta_{dg}(R,d)\cap\ker(d).$$
\end{Prop}

Proof.
Let us compare $\zeta(R)$ and $\zeta_{dg}(R,d)$. For $a\in \zeta(R)$ we need to have that
$\rann(a)$ is essential in $R$. For $a\in\zeta_{dg}(R,d)$ we need to have that $\rann(a)$
is a dg-essential dg-ideal. If $a\in\ker(d)$, then $\rann(a)$ is a dg-ideal, and
essential ideals are trivially dg-essential.
Hence $$\zeta(R)\cap\ker(d)\subseteq \zeta_{dg}(R,d)\cap\ker(d).$$
Further, if $I$ is a right ideal of $R$, then $I\cap\ker(d)$ is a right ideal of $\ker(d)$.
Hence, if $a\in\zeta(\ker(d))$, then $\rann_{\ker(d)}(a)$ intersects non trivially any
non zero ideal of $\ker(d)$. Likewise, if $a\in \zeta(\ker(d))$, then
$\rann_R(a)$ intersects non trivially all ideals $I$ of $R$ with $I\cap\ker(d)\neq 0$.
However, all dg-ideals $(I,d)$ of $(R,d)$ do intersect with $\ker(d)$.
Indeed, if $y\in I\neq 0$, then $d(y)\in I$ as well, and hence either $y\in\ker(d)$
or $d(y)\in\ker(d)$.  Hence
$$\zeta(\ker(d))\subseteq\zeta_{dg}(R,d)\cap\ker(d).$$
This shows the statement. \dickebox

\medskip

In Proposition~\ref{singularidealofkernel} we considered the (right) singular ideal of
the subalgebra $\ker(d)$ of the dg-ring $(R,d)$. Since we have a
surjective ring homomorphism $\ker(d)\lra H(R,d)$, we can also consider the singular ideal of
the ring $H(R,d)$, and the singular ideal  of the $\ker(d)$-module $H(R,d)$.
In the next lemma we are going to explore their relation.

\begin{Prop} \label{zetamoduleversusringunderpi}
Let $(R,d)$ be a differential graded ring. Let $\pi:\ker(d)\lra H(R,d)$ be the
natural homomorphism.
Then
\begin{enumerate}
\item
$\zeta(H(R,d))$ coincides with the singular submodule of the $\ker(d)$-module $H(R,d)$.
\item
$\pi(\zeta(\ker(d)))\subseteq \zeta(H(R,d))$.
\end{enumerate}
\end{Prop}

Proof. Since $\pi$ is surjective, $\pi$ induces a bijection between
the ideals of $H(R,d)$ and the ideals of $\ker(d)$ containing $\ker(\pi)=\im(d)$.
Hence,
\begin{eqnarray*}
h\in\zeta(H(R,d))&\Leftrightarrow&\rann_{H(R,d)}(h)\textup{ is essential in }H(R,d)\\
&\Leftrightarrow& \forall_{0\neq I\leq_r H(R,d)}\;:\;I\cap\rann_{H(R,d)}(h)\neq 0
\end{eqnarray*}
Moreover, the annihilator $\rann_{H(R,d)}(h)$ of $h$ in $H(R,d)$ coincides with
the image under $\pi$ of  the annihilator of $h$ as a $\ker(d)$-module. Since
$\ker(\pi)$ certainly annihilates $h$, we get that  $\zeta(H(R,d))$ coincides
with the singular submodule of the $\ker(d)$-module $H(R,d)$.

Let $a\in\zeta(\ker(d))$. This is equivalent to
$a\in \ker(d)$ and $\rann_{\ker(d)}(a)$ essential in $\ker(d)$. But if  an element $b$ in
$\ker(d)$ annihilates $a$, then $\pi(0)=\pi(ab)=\pi(a)\pi(b)$, and hence
$\pi(\rann_{\ker(d)}(a))\subseteq\rann_{H(R,d)}(\pi(a))$.
Therefore if $\rann_{\ker(d)}(a)$ is essential in $\ker(d)$, then
$\pi(\rann_{\ker(d)}(a))$ is essential in $H(R,d)$.
\dickebox

\begin{Rem}
Note that the proof of Proposition~\ref{zetamoduleversusringunderpi}
shows that whenever $\pi:R\lra S$ is a surjective ring homomorphism, then
$\zeta(S)$ coincides with the singular $R$-submodule of the $R$-module $S$.
\end{Rem}

\begin{Lemma} \label{zetadgisdgideal}
Let $(R,d)$ be a differential graded ring. Then
$\zeta_{dg}(R,d)$ is a differential graded 
left ideal of $(R,d)$.
\end{Lemma}

Proof. Let $a,b\in\zeta_{dg}(R,d)$. Then there are dg-essential dg-right
ideals $E_a$ and $E_b$ such that
$aE_a=0=bE_b$.
By Lemma~\ref{essentiallemma}.\ref{essentiallemmaitem1} also $E_a\cap E_b$
is an essential dg-right ideal of $R$. Then $a$ and $b$ annihilate $E_a\cap E_b$, and hence also
$a-b$. Let $x\in R$. Then $xa$ annihilates $E_a$ as well, and hence $\zeta_{dg}(R,d)$
is stable by left multiplication with elements in $R$.
Now, for any homogeneous $x\in E_a$, supposing
that $a\in\zeta_{dg}(R)$ is homogeneous,  we have
$$d(a)\cdot x=d(ax)-(-1)^{|a|} a\cdot d(x)$$
and since $x\in E_a$, we have $ax=0$, whence also $d(ax)=0$.
Since $E_a$ is a dg-ideal, also $d(x)\in E_a$ and hence $a\cdot d(x)=0$. Therefore $d(a)$
annihilates $E_a$ as well, and hence $d(a)\in\zeta_{dg}(R,d)$. This shows the lemma.
\dickebox

\begin{Rem}\label{zetadgrelativetozeta}
We cannot show in general that $\zeta_{dg}(R,d)\subseteq \zeta(R)$.
Indeed, if $a\in\zeta_{dg}(R,d)$, then there is a dg-essential dg-ideal $(E,d)$
of $(R,d)$ such that $aE=0$. However, a dg-essential dg-ideal does not need to
be essential.
Since an essential ideal does not need to be a dg-ideal, we cannot show
$\zeta_{dg}(R,d)\supseteq \zeta(R)$ neither.
\end{Rem}

\begin{Rem}
Since $\zeta_{dg}(R,d)$ is a dg-left ideal of $(R,d)$, it is tempting to consider its homology
$H(\zeta_{dg}(R,d),d)$.
An element in $H(\zeta_{dg}(R,d),d)$ is represented by $y\in\ker(d)$ such that
$\rann(y)$ is dg-essential in $(R,d)$.
Since by Proposition~\ref{singularidealofkernel} we have
$$\zeta(\ker(d))\subseteq\zeta_{dg}(R,d)\cap \ker(d)$$
and since $H(\zeta_{dg}(R,d),d)$ is a quotient of $\zeta_{dg}(R,d)\cap \ker(d)$,
there is a natural map
$$\zeta(\ker(d))\lra H(\zeta_{dg}(R,d),d)$$
induced by the natural map
$$\ker(d|_{\zeta_{dg}(R,d)})\lra H(\zeta_{dg}(R,d),d).$$
Likewise, since by Proposition~\ref{singularidealofkernel} we have
$$\zeta(R)\cap\ker(d)\subseteq\zeta_{dg}(R,d)\cap \ker(d)$$
we also get a natural map
$$\zeta(R)\cap\ker(d)\lra H(\zeta_{dg}(R,d),d).$$
However, since for an element $a\in\ker(d)$ the property for the
dg-right ideal $\rann(a)$ to be essential is a lot more restrictive than
to be dg-essential, there is no hope to have surjectivity of either one of these maps.
An example is given below in Example~\ref{exbelow}.
\end{Rem}

\begin{Example}\label{exbelow}
Indeed, recall Example~\ref{examplematrixring}.
For a field $K$ we defined a structure of a dg-algebra on $A=Mat_2(K)$.
Then, $\zeta(A)=0$ and
$$\ker(d)=K\cdot\left(\begin{array}{cc}0&1\\0&0\end{array}\right) +
K\cdot\left(\begin{array}{cc}1&0\\0&1\end{array}\right).$$
Moreover, $$\zeta_{dg}(A)=\left(\begin{array}{cc}0&K\\0&K\end{array}\right)\textup{
and }H(\zeta_{dg}(A,d),d)=K\cdot \left(\begin{array}{cc}0&0\\0&1\end{array}\right)\simeq K.$$
Therefore in this case the map
$$\zeta(A)\cap\ker(d)\lra H(\zeta_{dg}(A,d),d).$$
is not surjective.

Further, in this case the only non trivial ideal of $\ker(d)$ is
$J:=K\cdot\left(\begin{array}{cc}0&1\\0&0\end{array}\right) $.
Hence, this ideal is essential. Its right annihilator ideal is $J$ itself, and hence
$\zeta(\ker(d))=J$. This shows that the map $$\zeta(\ker(d))\lra H(\zeta_{dg}(R,d),d)$$
is the zero map, which is neither surjective nor injective.
\end{Example}

\section{dg-Goldie rings; left and right dg-annihilators}

\label{dgleftrightannihilsect}

\begin{Def}
A differential graded ring $(R,d)$ is called {\em dg-left Goldie} if $(R,d)$ satisfies
\begin{itemize}
\item
the ascending chain condition on dg-left annihilators, and
\item
has finite left dg-uniform dimension.
\end{itemize}
Analogously we define {\em dg-right Goldie} rings.
\end{Def}

The statements of \cite[Proposition 2.14]{McconnellRobson} are formal, except the last statement,
and can be transposed
to the differential graded situation. Let us check this.

\begin{Prop}\label{firstpropertiesofdgsemiprime}
Let $(R,d)$ be dg-semiprime and let $(I,d)$ be a two-sided dg-ideal. Suppose $R\neq I$.
Then
\begin{enumerate}
\item\label{firstpropertiesofdgsemiprime1}
We have
$\lann(I,d)=\rann(I,d)=:\ann(I,d)$ and this is a two-sided dg-ideal.
\item \label{firstpropertiesofdgsemiprime2}
$\ann(I,d)$ is a two-sided dg-ideal with $I\cap \ann(I,d)=0$, and it is the unique one
which is maximal with respect to this property. In particular, $I\oplus \ann(I,d)$ is dg-essential in $R$.
\item\label{firstpropertiesofdgsemiprime3}
If $(R,d)/\ann(I,d)$ is strongly dg-semiprime, then
$\ann(I,d)$ is the intersection of those minimal dg-prime ideals which do not contain $I$.
\item\label{firstpropertiesofdgsemiprime4}
If the $(R,d)-(R,d)$-bimodule $(I,d)$ is dg-uniform then $\ann(I,d)$ is a minimal 
dg-prime ideal.  If  $R/\ann(I)$ is strongly dg-semiprime, then the converse holds as well.
\item\label{firstpropertiesofdgsemiprime5}
$(I,d)$ is a dg-essential two-sided dg-ideal if and only if $\ann(I,d)=0$.
\item \label{firstpropertiesofdgsemiprime6}
Suppose that $R/\ann(I)$ is strongly dg-semiprime.
If $(I,d)$ is not contained in any minimal dg-prime ideal, then $\ann(I,d)=\Prad(I,d)$.
In particular, if $(R,d)$ is strongly dg-semprime, then $(I,d)$ is dg-essential 
if and only if $(I,d)$ is not contained in any dg-prime ideal.
\end{enumerate}
\end{Prop}

Proof. We shall prove (\ref{firstpropertiesofdgsemiprime1}).
$X:=\{r\in R\;|\;rI=0\}$ is a two-sided dg-ideal by Lemma~\ref{annihilatrosaredgideals}.
Now, $$(IX)^2=I(XI)X=0$$ and $IX$ is a nilpotent
dg-ideal. Since $(R,d)$ is dg-semiprime, we get
$\dgnil(R,d)=0$, and therefore $IX=0$ and $X\subseteq  \{r\in R\;|\;Ir=0\}$. The same
argument holds the other way round, and hence we get a double inclusion
$$\{r\in R\;|\;rI=0\}\subseteq \{r\in R\;|\;Ir=0\}\subseteq \{r\in R\;|\;rI=0\}.$$
This shows (\ref{firstpropertiesofdgsemiprime1}).

\medskip

 We shall prove (\ref{firstpropertiesofdgsemiprime2}).
Let $J$ be a dg-ideal with $I\cap J=0$. Then $0=I\cap J\supseteq IJ$ and hence $J\subseteq \ann(I)$.
As $(I\cap (\ann(I))^2=0$, and as $(R,d)$ is dg-semiprime, we have $\dgnil(R,d)=0$, whence
we have $I\cap \ann(I)=0$. This shows (\ref{firstpropertiesofdgsemiprime2}).
The fact that $I\oplus \ann(I)$ is dg-essential is a consequence of
Lemma~\ref{essentiallemma}.(\ref{essentiallemmaitem6}).

\medskip

 We shall prove (\ref{firstpropertiesofdgsemiprime3}).
If $(P,d)$ is a minimal dg-prime ideal which does not contain $I$, then
$$I\cdot \ann(I,d)=0\subseteq P.$$
Since $P$ is dg-prime, either $I\subseteq P$ or $\ann(I,d)\subseteq (P,d)$.
The first case was excluded, and so the intersection $D$ of all (minimal) dg-prime ideals
not containing $I$ contains $\ann(I,d)$. 
Since $R/\ann(I)$ is strongly dg-semiprime, the intersection of all dg-primes of $R/\ann(I)$
is $0$. Taking preimages under $R\lra R/\ann(I)$ we get that the intersection of all 
dg-primes containing $\ann(I)$ is $\ann(I)$. Since $I\cap\ann(I)=0$, we have that $D$ 
equals the intersection of all dg-primes containing $\ann(I)$. Hence
we get $I\cap D=0$. 

Therefore, $D= \ann(I,d)$ by (\ref{firstpropertiesofdgsemiprime2}), whence the statement of
(\ref{firstpropertiesofdgsemiprime3}).

\medskip

We shall prove item (\ref{firstpropertiesofdgsemiprime4}).
Suppose that $(I,d)$ is dg-uniform.
Let $S$ and $T$ be two-sided dg-ideals with $ST\subseteq \ann(I,d)$. Then
$IST=0=STI$. If $IS=0$, then $S\subseteq \ann_r(I,d)= \ann(I,d)$ and we are done.
Likewise, if $TI=0$, then $T\in \ann_\ell(I,d)=\ann(I,d)$, and we are done as well.
If $IS\neq 0\neq TI$, then
$T\subseteq \ann_r(IS)=\ann(IS)$ and $S\subseteq \ann_\ell(TI)=\ann(TI)$.
Further, $0\neq IS\subseteq I\cap S$ and $0\neq TI\subseteq T\cap I$.
By item (\ref{firstpropertiesofdgsemiprime2}) $\ann(J,d)$ is a two-sided ideal with
$J\cap \ann(J,d)=0$, and it is maximal with this property, and furthermore the unique maximal one.

Further, $IS\subseteq I$ and $TI\subseteq I$.
If $0\neq (J,d)\subseteq (I,d)$, then
$\ann(I,d)\subseteq \ann(J,d)$ by definition. We claim that  $\ann(I,d)= \ann(J,d)$.
Indeed, $(I,d)$ is dg-uniform, hence $(J,d)$ is a dg-essential submodule of $(I,d)$.
Now, if $\ann(I,d)\subsetneq \ann(J,d)$, then $I\cap \ann(J,d)\neq 0$, since else
this would contradict  item (\ref{firstpropertiesofdgsemiprime2}), namely
the maximality of $\ann(I,d)$ as being maximal with $I\cap \ann(I,d)=0$.
But then $I\cap \ann(J,d)$ is a non zero dg-submodule of $I$, and since $(J,d)$ is dg-essential,
$J\cap (I\cap \ann(J,d))\neq 0$. However, $\ann(J,d)$ satisfies $J\cap \ann(J,d)=0$ by
item (\ref{firstpropertiesofdgsemiprime2}). This contradiction shows that $\ann(I,d)=\ann(J,d)$.
We can now consider $J=IS$ and this then implies $T\subseteq \ann(IS)=\ann(I)$. Hence $\ann(I,d)$ is dg-prime. 

Suppose that $\ann(I,d)$ is 
dg-prime. We shall need to see that $(I,d)$ is dg-uniform.
Let $(J,d)$ be a two-sided dg-ideal in $(I,d)$. We shall need to see that $(J,d)$
is dg-essential in $(I,d)$. Let $(K,d)$ be a dg-ideal in $(I,d)$. Then $J\cdot K\subseteq J\cap K$.
If $J\cdot K=0$, then $J\cdot K\subseteq\ann(I,d)$. Since $\ann(I,d)$ is dg-prime, 
either $J\subseteq \ann(I,d)$ or $K\subseteq\ann(I,d)$. 
However, $J\subseteq I$ and $K\subseteq I$ implies $J=0$ or $K=0$ by item (\ref{firstpropertiesofdgsemiprime2}).
Now, $\ann(J,d)$ is an intersection of minimal
dg-prime ideals, by item (\ref{firstpropertiesofdgsemiprime3}).
By definition $\ann(I,d)\subseteq \ann(J,d)$
as $J\subseteq I$. But $\ann(I,d)$ is a minimal dg-prime, whence contributes
to the intersection of minimal dg-primes giving $\ann(J,d)$. Therefore
$\ann(J,d)\subseteq \ann(I,d)$, and hence actually $\ann(J,d)=\ann(I,d)$ since the
other inclusion was seen above.
If now $(I,d)$ contains a direct sum of two two-sided dg-ideals $(I_1,d)\oplus (I_2,d)\subseteq (I,d)$,
then $I_1\cdot I_2\subseteq I_1\cap I_2=0$, and hence
 $\ann(I_1,d)$ contains $(I_2,d)$. Therefore $\ann(I_1,d)\supsetneq \ann(I,d)$. Taking
 $J=I_1$ in the discussion above, we get $\ann(I_1,d)=\ann(I,d)$.
This contradiction shows that the dg-uniform dimension is $1$, and hence, by
Proposition~\ref{sumofuniformsubmods}, $(I,d)$ is dg-uniform.
We proved item (\ref{firstpropertiesofdgsemiprime4}).

\medskip

We shall prove item (\ref{firstpropertiesofdgsemiprime5}).
If $(I,d)$ is dg-essential, by item (\ref{firstpropertiesofdgsemiprime2})
we need to have $\ann(I,d)=0$. Let us prove the other direction. Suppose that $\ann(I,d)=0$.
Let $(J,d)$ be a two-sided dg-ideal of $(R,d)$. Then $I\cdot J\subseteq I\cap J$. If $(I,d)$
is not dg-essential, then there is a non zero two-sided dg-ideal $(J,d)$ with $I\cap J=0$. Hence, by
item (\ref{firstpropertiesofdgsemiprime2}) we have that $J\subseteq \ann(I,d)$. But $\ann(I,d)=0$,
and this contradiction gives  item (\ref{firstpropertiesofdgsemiprime5}).

\medskip

We shall prove item (\ref{firstpropertiesofdgsemiprime6}).
By item (\ref{firstpropertiesofdgsemiprime3}) and the hypothesis,
we get that $\ann(I,d)$ is the intersection
of all dg-prime ideals, which is $\Prad(R,d)$.
In case $(R,d)$ is strongly dg-semiprime, then $\Prad(R,d)=0$. Hence
$(I,d)$ is not contained in any dg-prime ideal if and only if
$\ann(I,d)=0$ if and only if $(I,d)$ is dg-essential.
\dickebox

\medskip

\begin{Rem}
Recall that in the non dg-situation the intersection of all prime ideals
is the nil radical, which is
$0$ for semiprime rings. Hence $I$ is essential in the situation of
item (\ref{firstpropertiesofdgsemiprime6}). In the dg-situation we do not have this
property (cf Example~\ref{dgPraduneaqualdgnil}) and we need the stronger hypothesis
of strongly dg-semiprime rings.
\end{Rem}

Let $(R,d)$ be a dg ring. Let $(I,d)$ be a two-sided dg-ideal.
Trivially, if $I$ is essential, then $(I,d)$ is also dg-essential.
The converse sometimes holds as well in a very specific situation.

\begin{Cor}
If $(R,d)$ is a differential graded algebra, and suppose that $R$ is semiprime as a ring.
Then a two-sided dg-ideal $(I,d)$ is dg-essential if and only if $\ann(I,d)=0$, if and only if
$\ann(I)=0$, if and only if $I$ is essential.
\end{Cor}

Proof.
$R$ is assumed to be semiprime, and hence $R$ does not contain a non zero two-sided
nilpotent ideal, whence neither a two-sided nilpotent dg-ideal. Hence $(R,d)$ is dg-semiprime as well.
We may now apply Proposition~\ref{firstpropertiesofdgsemiprime}
item (\ref{firstpropertiesofdgsemiprime5})
and its classical counterpart \cite[Proposition 2.14 item (5)]{McconnellRobson}.
If $(I,d)$ is dg-essential, then $\ann(I,d)=0$ and since by definition $\ann(I,d)=\ann(I)$,
we have that \cite[Proposition 2.14 item (5)]{McconnellRobson} implies that $I$ is essential. \dickebox

\medskip

As in the classical situation \cite[2.2.3 and 2.2.10]{McconnellRobson} we may prove

\begin{Lemma}
Let $(A,d)$ be a differential graded ring and let $(M,\delta)$ be a differential graded
$(A,d)$-module.
\begin{itemize}
\item
If $(N,\delta)$ is a dg-complement dg-submodule of $(M,\delta)$, then for
all dg-submodules $(L,\delta)$ of $(M,\delta)$  
with $N\subsetneq L$ there is a non zero dg-submodule $(S,\delta)$
of $(L,\delta)$ with $S\cap N=0$.
\item
Then the following are equivalent:
\begin{itemize}
\item $\dgudim(M,\delta)<\infty$
\item $(M,\delta)$ satisfies the ascending
chain condition on dg-complement submodules. The dg-uniform dimension of $(M,\delta)$ is the
maximal length of an ascending chain of dg-complement dg-submodules.
\end{itemize}
\end{itemize}
\end{Lemma}

Proof. As for the first item, suppose that $(X,\delta)$ is a dg-submodule of $M$, such that
$N$ is a dg-complement to $X$. Then put $S:=L\cap X$. Then $(S,\delta)$ is a dg-submodule
of $(M,\delta)$ since $L$ and $X$ are  dg-submodules
of $(M,\delta)$. Further,
$$ S\cap N=N\cap L\cap X=N\cap X=0 $$
since $N\subseteq L$ and since $X$ is a dg-complement to $N$. If $S=0$, then $L$ would be
the dg-complement to $X$, which is excluded since by hypothesis the dg-complement $N$ of $X$ is
strictly smaller than $L$.

\medskip

As for the second item, let $n=\dgudim(M,\delta)$. Let
$$0<S_1<S_2<S_3<\dots<S_t$$
be a maximal chain of dg-complement dg-submodules of $(M,\delta)$.
We claim that $t\leq n$. Indeed, applying the statement of the 
first item to $S_1<S_2$, we obtain a dg-submodule $S_2'$
of $S_2$ with $S_1+S_2=S_1\oplus S_2'$. Similarly, $S_2<S_3$ gives a dg-submodule $S_3'<S_3$ such
that $S_2+S_3'=S_2\oplus S_3'$. Since $S_2$ contains $S_1\oplus S_2'$, we obtain a direct sum
$S_1\oplus S_2'\oplus S_3'$ of three dg-submodules. By induction we get a direct sum of $t$
non zero dg-submodules $S_i'\leq S_i$ for all $i\in\{1,\dots,t\}$ such that
$S_1\oplus S_2'\oplus\dots\oplus S_t'$ is a direct sum of dg-submodules of $(M,\delta)$.
Since $\dgudim(M,\delta)=n$, we get $t\leq n$.

If we have a direct sum $\bigoplus_{i=1}^tM_i$ of dg-submodules of $(M,\delta)$,
then, for $s<t$ the dg-complement of  $\bigoplus_{i=1}^sM_i$ contains  $\bigoplus_{i=s+1}^tM_i$
but does not contain $M_i$ for any $i\leq s$.
This direct sum hence induces a chain of dg-complement dg-submodules of length at least $t$.
This then proves the second item. \dickebox

\medskip

Recall from Theorem~\ref{orelocalisationofdg} that a differential graded structure
on a ring $R$ can be extended to the Ore localisation at a set of homogeneous and
regular Ore set. Moreover, since the Ore set is formed by regular elements,
the natural homomorphism to the localisation is injective.
Compare the following lemma with \cite[2.12]{McconnellRobson}.
We identify $R$ with its image in $R_S$ under the canonical dg-homomorphism
$R\stackrel{\lambda}\lra R_S$.

\begin{Lemma}\label{dgudimforlocalisationlemma}
Let $(R,d)$ be a differential graded ring and let $S$ be an Ore set
of homogeneous regular elements of $R$.  Let $(I,d)$ be a dg-right ideal of $(R,d)$, and let
$(J,d)$ be a dg-right ideal of $(R_S,d)$.
\begin{enumerate}
\item \label{item1local} Then $(I,d)$ is dg-essential in $(R,d)$ if and only if $(I\cdot R_S,d_S)$ is dg-essential in $(R_S,d_S)$.
\item \label{item2local} Then $(J,d_S)$ is dg-essential in $(R_S,d_S)$ if and only if $(J\cap R,d)$ is dg-essential in $(R,d)$.
\item \label{item3local} $\dgudim_{(R,d)}(I)=\dgudim_{(R_S,d_S)}(I\cdot R_S)$
\item \label{item4local} $\dgudim_{(R,d)}(J\cap R)=\dgudim_{(R_S,d_S)}(J)$
\end{enumerate}
\end{Lemma}

Proof.
We shall prove item (\ref{item1local}).
Suppose that $(I,d)$ is dg-essential in $(R,d)$. Let $(L(S),d_S)$ be a dg-ideal in
$R_S$. Then $(L(S)\cap R,d)$ is a dg-ideal in $(R,d)$ and hence $L(S)\cap R\cap I=L(S)\cap I\neq
0$. Then $(L(S)\cap I)\cdot R_S\subseteq L(S)\cap I\cdot R_S$, and since $S$ is formed
by regular elements, $\lambda$ is injective, whence $L(S)\cap I\cdot R_S\neq 0.$
This shows that $I\cdot R_S$ is dg-essential.

Suppose now that $I\cdot R_S$ is dg-essential and let $L$ be a dg-ideal of $R$.
If $I\cap L=0$, then $I+L=I\oplus L$ is a dg-right ideal of $R$, and hence, using
flatness of the localisation,
$$(I+L)R_S=(I\oplus L)R_S=(I\cdot R_S)\oplus (L\cdot R_S)$$
is a dg-right ideal of $R_S$. Since $L\cdot R_S$ is a dg-ideal of $R_S$, and since the intersection
with $I\cdot R_S$ in $0$, this contradicts the fact that $I\cdot R_S$ is dg-essential.

\medskip

We shall prove item (\ref{item2local}).
Suppose that $(J,d_S)$ is dg-essential in $(R_S,d_S)$ and let $(L,d)$ be a
dg-ideal in $(R,d)$. Then $J\cap R$ is a dg-ideal in $R$ and if $(J\cap R)\cap L=0$, then
since $(J\cap R)\cap L=J\cap L$, we also get $J\cap (L\cdot R_S)=0$. This implies $L\cdot R_S=0$
since $L\cdot R_S$ is a dg-ideal of $R_S$ and $J$ is dg-essential. However,
$L\cdot R_S=0$ implies $L=0$. Hence $J\cap R$ is dg-essential.

Suppose that $J\cap R$ is dg-essential in $R$ and let $L$ be a dg-ideal of $R_S$. If
$J\cap L=0$, then $J\cap L\cap R=(J\cap R)\cap (L\cap R)$. Since $J\cap R$ is
dg-essential in $R$, we get $L\cap R=0$. This implies $L=0$ and we showed
that $J$ is dg-essential.


\medskip

We shall prove item (\ref{item3local}).
We have an additive functor
\begin{eqnarray*}
(R,d)-dg-mod&\lra&(R_S,d_S)-dg-mod\\
M&\mapsto&M\otimes_R(R_S,d_S)
\end{eqnarray*}
This functor preserves direct sums, and hence
$\dgudim_{(R,d)}(I)\leq \dgudim_{(R_S,d_S)}(I\cdot R_S)$.
If $(I,d)$ is a dg-uniform ideal of $(R,d)$, then $(I\cdot R_S,d_S)$ is dg-uniform. This follows from 
\cite[(1.16) Proposition]{McconnellRobson}.
Using item (\ref{item1local}) we have that a direct sum $\bigoplus I_i$ of uniform dg-ideals of $(R,d)$ 
is dg-essential if and only if the direct sum $\bigoplus I_i\cdot R_S$ is dg-essential.
Hence we proved the statement.

\medskip

We shall prove item (\ref{item4local}).
Let $\bigoplus_{i=1}^n J_i\subseteq J$ be a direct sum of dg-uniform ideals of $R_S$.
Then the dg-$R$-ideal $\sum_{i=1}^n (J_i\cap R)$ in $J\cap R$ is actually a direct sum, and hence
$\dgudim_{(R,d)}(J\cap R)\geq \dgudim_{(R_S,d_S)}(J)$. Using \cite[(1.16) Proposition]{McconnellRobson}
again we see that if $J_i$ is dg-uniform, then also $J_i\cap R$ is dg-uniform. 
Further, by item (\ref{item2local}) we have that $\bigoplus_{i=1}^n J_i$  is
dg-essential in $J$ if and only if $\sum_{i=1}^n (J_i\cap R)$ is essential in $J\cap R$.
Hence we proved the statement. \dickebox


\begin{Rem}
The fact that \cite[(2.12) Lemma]{McconnellRobson} follows from
\cite[(1.16) Proposition]{McconnellRobson}, and that theses statements are independent of the
presence of a dg-structure, also the generalisation  to the differential graded situation
Lemma~\ref{dgudimforlocalisationlemma} follows from \cite[(1.16) Proposition]{McconnellRobson}.
\end{Rem}

\begin{Rem}\label{maximalannihilators}
If $I$ is a dg-ideal such that $\lann(I)$ is maximal within the set of left annihilators,
then $I=Ra$ for some $a\in\ker(d)$. Indeed, let $0\neq a\in I$, and if $d(a)\neq 0$,
then replace $a$ by $d(a)$, which is again in $I$ since $I$ is a dg-ideal.
Hence, we may find $a\in I\cap\ker(d)$.
Then $Ra$ is a dg-ideal and $Ra\subseteq I$. This implies $\lann(I)\subseteq \lann(Ra)$.
Maximality of $\lann(I)$ shows that left annihilators which are maximal within the set of
let annihilators are annihilators of dg-principal ideals.
\end{Rem}

Recall that in the classical ungraded case and differential $0$ we have the following lemma.

\begin{Lemma} \cite[Lemma 3.2]{McconnellRobson}
Let $R$ be a ring and suppose that $R$ satisfies the
ascending chain condition on left annihilators. Then
\begin{enumerate}
\item Each maximal left annihilator has the form $\lann(a)$ for some $a\in R$.
\item For any $b\in R$ there is an integer $m$ such that $\lann(b^n)=\lann(b^m)$
for any $n\geq m$. Then, for these $n\geq m$ we have $\lann(b^n)\cap Rb^n=0$
\item Each non zero nil left ideal contains a non zero nilpotent ideal.
\end{enumerate}
\end{Lemma}

We consider the differential graded case.

\begin{Cor}
Let $(R,d)$ be a differential graded ring, and suppose that $R$ satisfies the
ascending chain condition on left annihilators. Let $b\in R$. Then
If  $\lann(b^n)=\lann(b^m)$ for any $n\geq m$, and if $b^n\in\ker(d)$,
then $Rb^n\oplus \lann(b^n)$ is an essential dg-ideal, whence in particular dg-essential.
\end{Cor}

Proof. By \cite[Lemma 3.3]{McconnellRobson} we have that $Rb^n\oplus \lann(b^n)$ is an
essential left ideal. Since $b^n\in\ker(d)$, $Rb^n$ is a dg-left ideal,
and $\lann(b^n)$ is a dg-left ideal by Lemma~\ref{annihilatrosaredgideals}.
This shows the statement. \dickebox

\medskip

We are generalising \cite[Lemma 3.4]{McconnellRobson} to the dg-situation.

\begin{Lemma}
Let $(R,d)$ be a dg-ring
\begin{itemize}
\item
If $c\in \ker(d)$ is right regular,
then the dg-right ideal $(cR,d)$ is dg-essential in $(R,d)$
\item
Suppose that $(R,d)$ is  dg-semiprime, suppose that it has finite left dg-uniform
dimension, and suppose that $\zeta_{dg}(R,d)=0$. Then if $c\in \ker(d)$ is right regular,
$c$ is regular.
\end{itemize}
\end{Lemma}

Proof.
Since $c\in\ker(d)$ is right regular, we get $cx=0$ implies $x=0$, and also $cR$ is a dg-ideal.
Hence
\begin{align*}
R\ra&R\\
x\mapsto&cx
\end{align*}
is an isomorphism of $(R,d)$-right dg-modules.
Hence, using the differential graded uniform dimension of right ideals, we get
$\dgudim(R)=\dgudim(cR)$, and therefore using Lemma~\ref{dgudimofsubmodules} item ii)
(transposed to right modules)
we get $(cR,d)$ is dg-essential in $(R,d)$.

We need to see that $c$ is left regular.
Since $\zeta_{dg}(R,d)=0$, there is no non zero element $c$ and a dg-right ideal $E$ such that
$cE=0$. Since $cR$ is dg-essential, we need to have $\lann(cR)=0$. Since
$\lann(c)\subseteq \lann(cR)$, we also obtain $\lann(c)=0$. But this shows that
there is no non zero element $y\in R$ such that $yc=0$. This is tantamount to say
that $c$ is left regular. \dickebox

\section{Differential graded Goldie-theorem}

\label{dgGoldietheoremsect}

Recall that a dg-algebra $(A,d)$ is dg-prime if whenever $(I,d)$ and $(J,d)$ are two-sided
dg-ideals of $(A,d)$, then $IJ=0$ implies $I=0$ or $J=0$. A graded prime (or gr-prime) algebra
$A$ is a dg-algebra $(A,d)$ which is dg-prime for $d=0$.

\begin{Lemma}\label{grprimeverssusdgprime}
Let $(A,d)$ be a dg-algebra. Then if $\ker(d)$ is gr-prime, we have that $(A,d)$  is dg-prime.
\end{Lemma}

Proof. Let $(A,d)$ be a dg-algebra and let $S:=\ker(d)$.
Let $(I,d)$ and $(J,d)$ be two two-sided dg-ideals. Then $IJ=0$ implies
$(I\cap S)\cdot(J\cap S)=0$ in $S$, and since $S$ is gr-prime, we get $I\cap S=0$ or $J\cap S=0$.
However, for any $x\in I$ we get if $x\not\in\ker(d)$, then $0\neq d(x)\in\ker(d)\cap I$.
Likewise for $J$.
Hence $(I,d)=0$ or $(J,d)=0$. This shows the statement. \dickebox

\begin{Lemma}\label{grunifordimversusdguniformdim}
Let $(A,d)$ be a differential graded algebra. If $\ker(d)$ has finite gr-uniform dimension, then
$(A,d)$ has finite dg-uniform dimension.
\end{Lemma}

Proof. Put $S:=\ker(d)$.
 If $I_1\oplus I_2\oplus\dots\oplus I_n$ is a direct sum of
two-sided dg-ideals of $(A,d)$, then  $(I_1\cap S)\oplus (I_2\cap S)\oplus\dots\oplus (I_n\cap S)$
is a direct sum of two-sided graded ideals in $S$. Again, as in Lemma~\ref{grprimeverssusdgprime}
we see that $I\cap S=0$ implies $I=0$. This shows the lemma. \dickebox

\begin{Cor}
Let $(A,d)$ be a differential graded algebra and suppose that $(A,d)$ is two-sided
dg-Noetherian (i.e. all increasing sequences of two-sided dg-ideals becomes stationnary).

If $\ker(d)$ is left gr-Goldie, then $(A,d)$ is left dg-Goldie.
\end{Cor}

Proof. Lemma~\ref{grunifordimversusdguniformdim} and the hypothesis on the Noetherianity
shows the statement. \dickebox

\begin{Theorem} \label{essentialcontainsregular}
Let $R$ be a commutative ring and let $(A,d)$ be a differential graded $R$-algebra.
Suppose that $\ker(d)$ is a gr-prime ring and suppose
that $\ker(d)$ is right gr-Goldie. {\red Suppose that 
\begin{enumerate}
\item either the homogeneous regular elements $S_A$
of $A$ form an Ore set of $A$ 
\item or else the homogeneous regular elements $S_{\ker(d)}$ of $\ker(d)$ form an Ore set in $A$.
\end{enumerate}}

{\red Then $S_{\ker(d)}\subseteq S_A$ and hence the natural homomorphism from $A$ 
to the localisation at $S_A$ or at $S_{\ker(d)}$ is injective. }
Further, {\red in case (1) } the localisation of $(A,d)$ at {\red $S_A$ and in case (2) the 
localisation of $(A,d)$ at  $S_{\ker(d)}$}  is a
dg-simple differential graded $R$-algebra (in the sense that
there is no non zero non trivial two-sided dg-ideal).
\end{Theorem}

Proof. Let $(I,d)$ be a differential graded two-sided ideal of $(A,d)$. 
{\red Then by Lemma~\ref{essentiallemma} and the fact that $(A,d)$ is dg-prime
we get that $(I,d)$ is dg-essential.}
Hence
by Theorem~\ref{GoodearStafford} we know that $\ol I:=I\cap \ker(d)$
is either $0$ or, {\red since $\ker(d)$ is gr-prime, $\ol I$ is gr-essential, and} 
contains a homogeneous regular element.
$\ol I$ cannot be $0$ since for any $x\in I\setminus\{0\}$
we either have $x\in\ker(d)$, or else $d(x)\in I\cap\ker(d)$.
Hence $\ol I\subseteq I$ contains a homogeneous regular element $y$.
We claim that $y$ is regular in $A$ as well. Indeed, if $xy=0$ for some
$x\in A\setminus\ker(d)$. Considering each homogeneous component of $x$
separately, using that $y$ is homogeneous, we may assume that $x$
is homogeneous as well. Given hence a homogeneous $x\in A$ with $xy=0$, then
$$0=d(x\cdot y)=d(x)\cdot y+(-1)^{|x|} x\cdot d(y)=d(x)\cdot y$$
since $y\in\ker(d)$. But then $d(x)\in\ker(d)$, and $y$ regular in $\ker(d)$
implies that $d(x)=0$. Since $x\not\in\ker(d)$ by hypothesis we reach a contradiction.
Likewise, $y$ is right regular as well.
Using Theorem~\ref{orelocalisationofdg} we may localise at
the set of homogeneous regular elements $S$ and then $A_S\cdot I=A_S$.
Let now $(L,d)$ be a two-sided dg-ideal of $(A_S,d)$. Then $I:=L\cap A$
is a two-sided dg-ideal of $(A,d)$. By the above it contains a regular element,
and hence $L\supseteq A_S\cdot I=A_S$. Therefore $(A_S,d)$ does not contain any
proper non zero two-sided dg-ideal. We showed that $(A_S,d)$ is dg-simple.
\dickebox

\bigskip

We close with a lemma which is an analogue of the 'lying over property' in commutative algebra.

A graded ring is graded hereditary if any graded ideal is projective.

\begin{Lemma}\label{gradedhered}
Let $(A,d)$ be a differential graded algebra and let $S:=\ker(d)$.
Suppose that $S$ is graded hereditary. Then for any graded ideal $I$ of $S$
we have that $A\cdot I=:J$ is a differential graded ideal of $(A,d)$, and
$S\cap J=I$.
\end{Lemma}

Proof. Since $A$, $S$ and $I$ are graded, and since $A$ is a graded $S$-left module,
it is clear that $A\cdot I=:J$ is a graded ideal of $A$. Further, for any
homogeneous $a\in A$ and $y\in I\subseteq\ker(d)$, we get
$$d(a\cdot y)=d(a)\cdot y+(-1)^{|a|}a\cdot d(y)=d(a)\cdot y.$$
Hence $J$ is a dg-ideal of $(A,d)$.
Since $S$ is supposed to be graded hereditary, $I$ is projective, and
we may suppose that $Y:=\{y_i\;|\;i\in F\}$ is an $S$-basis of $I\oplus X$ for
some graded $S$-module $X$.

Recall that the differential of $A\otimes_S(I\oplus X)$ is $d\otimes_S\textup{id}_{I\oplus X}$.
Then let $y=\sum_{i=1}^na_i\otimes_S y_i\in J\cap S$ for elements $a_i\in A$.
Hence $$0=d(y)=\sum_{i=1}^nd(a_i)\otimes y_i$$
and since the set $Y$ is $S$-free, we infer $d(a_i)\in\ker(d)$.
Therefore we get that $d(a_i)=0$ for every
$i\in\{1,\dots,n\}$ and hence $a_i\in S$. This shows that $y\in I$.
\dickebox

\end{document}